\documentclass[11pt]{article}
\usepackage[utf8]{inputenc}
\usepackage{amsmath, amssymb, latexsym,amsfonts}

\usepackage{hyperref}

\usepackage{fancyhdr}
\numberwithin{equation}{section}
\newtheorem{theorem}{Theorem}[section]
\newtheorem{lemma}[theorem]{Lemma}

\title{Sums Associated with Orbits in the Binary Dynamical System}
\author{Rodney Nillsen \\
 University of Wollongong\\
 NSW 2522, Australia}
 \date{June 2023}
\begin{document}
\maketitle

\begin{center}
{\bf Abstract}
\end{center}
 \noindent{ In 1930,  G. H. Hardy and J. E. Littlewood  derived  results concerning  rates of divergence of certain series involving cosecants. In more recent terminology, one of their results can be interpreted in terms of  the behaviour of orbits in  a dynamical system that is a rotation on the unit circle.     Now, the expansion of numbers in $[0,1)$ to the base $2$ can be associated  with a different dynamical system -- the binary  system. This article considers  orbit behaviour in the binary system that corresponds  to the behaviour that was, in effect, observed by Hardy and Littlewood in  systems involving rotations.  Given a typical number in $[0,1)$,  the sequence of its binary digits  may be arranged as an infinite sequence of consecutive, non-empty, finite blocks, each block consisting  of all zeros  or all ones. The relationships  between the lengths of these blocks determine Hardy-Littlewood types of behaviour  associated with the number. Amongst other results, upper and lower estimates are derived for the sums of powers of the reciprocals of points in the orbit of the number.  These estimates are  in terms of the lengths of the associated blocks. A necessary and sufficient condition is found for the essential `equivalence' of the upper and lower estimates.  Almost all numbers in $[0,1)$ satisfy this condition. }


\section{Introduction}
 
 G. H. Hardy and J. E. Littlewood  \cite{HL}  considered the behaviour of some series involving  cosecants. Among the results they found was the following: \emph{there are irrational numbers  $\theta$  for which there is $K_{\theta}>0$  such that }

  \begin{equation}
  \sum_{k=1}^{n}\frac{1}{\sin^2{k\theta\pi}}\le K_{\theta}n^2, \ {\rm for\ \ all}\ n\in {\mathbb N},\label{eq:HL}
  \end{equation}
  where $\mathbb N$  is the set of natural numbers
Concerning this result they say: `This lemma is not actually used, and we include it because it is interesting in itself' \cite[p. 259]{HL}.

 This   result    can be expressed as a statement about orbits in the dynamical system given by rotations on the unit  circle.  Let's make the definition that a \emph{dynamical system} is a pair $(S,g)$, where $S$ is a set and $g$ is a function with $g:S\longrightarrow  S$. In this case, the composition of $g$ with itself, taken $k$ times, is denoted by $g^k$.  We take $g^0$ to be given by $g^0(x)=x$ for all $x$. If $x\in S$, the sequence $(x,g(x), g^2(x), \ldots)$ is called the \emph{orbit of $x$ in $S$ under $g$}, or simply the \emph{orbit} of $x$. Now, consider  the dynamical system $({\mathbb T}, \rho_{\theta})$, where ${\mathbb T}=\{z:z\in {\mathbb C} \ {\rm and }\ |z|=1\}$, $\theta\in [0,2\pi)$   and $\rho_{\theta}:{\mathbb T}\longrightarrow {\mathbb T}$ is given by 
$\rho_{\theta}(z)=e^{i\pi\theta}z.$ In this case, $\rho_{\theta}^k(z)=e^{ik\theta \pi}z$ for all $k\in \{0\}\cup {\mathbb N}$. Observe   that
$|1-\rho_{\theta}^k(1)|^2=|1-e^{ik\theta \pi}|^2\le 4,$
so that $ |1-\rho_{\theta}^k(1)|^{-2}\ge 1/4$, and so   the series 
$\sum_{k=1}^{\infty}|1-\rho_{\theta}^k(1)^{-2}|$ 
 is divergent.   As $|1-\rho_{\theta}^k(1)|^{-2}=(\sin^{-2} k\theta \pi/ 2)/4$, we can see 
  the  Hardy-Littlewood result (\ref{eq:HL}) as telling us that for some values $\theta\in (0,1)$ the rate of divergence of the series $\sum_{k=1}^{\infty}|1-\rho_{\theta}^k(1)|^{-2}$ is limited in the sense that there is $K_{\theta}>0$ such that
\begin{equation}
\sum_{k=1}^n\frac{1}{|1-\rho_{\theta}^k(1)|^2}\le K_{\theta}n^2, \ {\rm for \ all}\  n\in {\mathbb N}.\label{eq:HLproperty}
\end{equation}
 Here, we investigate  this interpretation of the Hardy-Littlewood result   in a different  dynamical system, the binary (or dyadic) system.

The binary  system is the dynamical system $([0,1), f)$, where 
\begin{equation}
f(x)=\begin{cases} 2x,\  {\rm if}\  0\le x<1/2, \ {\rm and} & \\ 2x-1,\ {\rm if}\  1/2\le x<1.\end{cases}\label{eq:dynamical}
\end{equation}
 If $x\in [0,1)$, $x$ is called a \emph{binary rational} if $x=0$ or $x=k/2^n$ for some $k,n\in {\mathbb N}$ with $1\le k< 2^n$. Motivated by the Hardy-Littlewood result  as  in  (\ref{eq:HLproperty}), for each number in $x\in [0,1)$ that is not a binary rational, and given $p>0$,  we will  estimate,  from   above and below,  the sum 
$\sum_{k=1}^n1/f^{k-1}(x)^p.$  
Functions $\Phi_1$ and $\Psi_1$ mapping $\mathbb N$ into $(0,\infty)$ are defined in terms of $p$ and  the lengths of the blocks or `runs' of zeros and ones in the binary  expansion of $x$,  and they have the property that there are $c_1,c_2>0$ such that
\begin{equation}
c_1\Phi_1(n)\le  \sum_{k=1}^n\frac{1}{f^{ k-1}(x)^p}\le c_2\Psi_1(n), \ {\rm for \ all}\ n\in {\mathbb N}.\label{eq:basicinequality}
\end{equation}
Also, corresponding functions $\Phi_2$ and $\Psi_2$ are defined having the property that for some $d_1,d_2>0$,
\begin{equation}
d_1\Phi_2(n)\le  \frac{1}{n}\sum_{k=1}^n\frac{1}{f^{ k-1}(x)^p}\le d_2\Psi_2(n), \ {\rm for \ all}\ n\in {\mathbb N}.\label{eq:basicinequality2}
\end{equation}
Relationships between these  estimates are discussed. A necessary and sufficient condition is  given in terms of the lengths of the blocks in the binary expansion of $x$  which ensures  the estimates in   (\ref{eq:basicinequality})  are `equivalent',  `optimal' or `sharp', in the sense that  $\Phi_1$ can replace $\Psi_1$ in the right hand side of (\ref{eq:basicinequality}). Also,  a corresponding necessary and sufficient condition is  given which ensures that the estimates in     (\ref{eq:basicinequality2}) are `sharp' in  the sense that $\Phi_2$ can replace $\Psi_2$. 

The condition that ensures that the estimates in (\ref{eq:basicinequality}) are `sharp' in the above sense imposes a restriction on the `regularity' of the digits in the binary expansion of $x$, as expressed in (\ref{eq:firstbound}) below. Also, the condition that ensures that the estimates in (\ref{eq:basicinequality2}) are `sharp' imposes a stronger, symmetric restriction on the `regularity' of the digits in the binary expansion of $x$, as expressed by (\ref{eq:regularbound}) and Theorem \ref{theorem2}. Now, the definitions of normal and simply normal numbers also involve notions of regularity in the digits of the number, and it turns out that both normal and simply normal numbers satisfy condition (\ref{eq:regularbound}). Thus, some results in this paper can be regarded as concerning an extension of the notions of normal and simply normal numbers, but where there is a weaker notion of `regularity' in the distribution of the digits. Whereas, if $x\in [0,1)$, the normality of $x$ is equivalent to the uniform distribution of  the points in the orbit sequence $(f^{n-1}(x))$ in $[0,1)$ (\cite[p. 70]{Ku}), the weaker notion of `regularity' in this paper is related to the  sharpness of the upper and lower bounds in the sums $\sum_{k=1}^n1/f^{k-1}(x)^p$, for some given index $p>0$. In the one case, `regularity'  consists in the uniform distribution of orbit points, while in the other case it consists in the coincidence of  the upper and lower bounding functions  of sums involving \emph{reciprocals} of the orbit points.

A particular problem  relating to the above concerns 
   the  numbers $x\in [0,1)$  that  have  the  property that  there is $K>0$, depending upon $x$,  such that 
\begin{equation}  
n\le \sum_{k=1}^n\frac{1}{f^{k-1}(x)^p}\le Kn,  {\rm for\  all}\  n\in {\mathbb N}. \label{eq:asharpcase}
\end{equation}
Note that  as $1\le 1/f^{k-1}(x) $, the left hand inequality in (\ref{eq:asharpcase}) holds for all  $n$.   It is shown that if $p>1$ and $x$ is a  simply normal number to base $2$,  the estimates in (\ref{eq:basicinequality}) and  (\ref{eq:basicinequality2})   are sharp,  in the sense mentioned above. In the case of these numbers the functions $\Phi_1$ and $\Phi_2$  take a simplified form. No normal number has  property  (\ref{eq:asharpcase}), but the set of  simply normal numbers that have property (\ref{eq:asharpcase}) is an uncountable set of measure zero. 

 \section{Preliminaries}
 We define $\Sigma$ to be the subset of  $[0,1)$ consisting of the numbers  that are not binary rationals.   Each $x\in[0,1)$  has a binary  expansion  to the base $2$:  for each $k\in {\mathbb N}$ there is  $d_k(x)\in \{0,1\}$ such that  
\[x=\sum_{k=1}^{\infty}\frac{d_k(x)}{2^k}.\]
When $x\in \Sigma$ the digits $d_k(x)$ are uniquely determined. So,  when $x\in \Sigma$,  we may refer to \emph{the} binary expansion of $x$ and   to \emph{the}   binary digits of $x$. In this case, the sequence of binary digits of $x$ contains no infinite sequence of consecutive zeros, and no infinite sequence of consecutive ones, and this property characterises the elements of $\Sigma$.  Details concerning the  expansion of numbers to a base may be found in \cite[pp. 92-101]{Ni} and \cite[pp. 64-67]{St}.  

The  dynamical system $([0,1), f)$ as in (\ref{eq:dynamical}) is  connected to the expansion of numbers to the base $2$.  This  seems  to have been realised first  by D. D. Wall \cite{Wa}.    One aspect of this connection is  that for all $x\in \Sigma$ and all $k\in {\mathbb N}$,
\[d_k(f(x))=d_{k+1}(x).  \]
That is,  $f$ shifts the sequence of binary digits of $x$ one position to the left, to give the sequence of digits of $f(x)$. We see that $f:\Sigma\longrightarrow \Sigma$ and that $d_j(f^{k}(x))=d_{j+k}(x)$. If one observes that $f(x)=0$ if and only if $x=0$ or $x=1/2$,  and that $f(x)$ is a binary rational if and only if $x$ is a binary rational, it can been seen that
\[\Sigma=\Bigl\{x:x\in [0,1)\ {\rm and}\ f^n(x)\ne 0, \ {\rm for \ all}\ n=0,1,2,\ldots\Bigr\}.\]
A   \emph{block} is any  finite sequence of zeros and ones,  and the \emph{length} of a block is its number of terms. A block may be empty, in which case its   length is zero. 

Let $x\in \Sigma$ be given. Let the   binary expansion  of $x$ be
 $x=\sum_{k=1}^{\infty}d_k/2^k.$
 The digit $d_k$ is said to have \emph{position} $k$.
 The sequence of digits $d_1,d_2,\ldots$ of $x$  is conventionally written as $d_1d_2d_3\ldots$ and  may be written as a sequence of juxtaposed blocks of zeros and ones    
 \begin{equation} C_0B_1C_1B_2C_2 \cdots,\label{eq:blocks}
 \end{equation}
 where each of $B_1,B_2,B_3,\ldots $ is a non-void finite sequence whose terms  are all $0$ and each of $C_0, C_1,C_2,C_3,\ldots,$ is a  finite sequence whose terms  are all $1$.  The sequence $C_0B_1C_1B_2C_2 \cdots$  in (\ref{eq:blocks}) is called the \emph{block decomposition} of $x$. 
    The \emph{length} of a block  is the number of terms  in the block. The block $C_0$ may be empty,  but all other blocks are non-empty. The length of the block $B_j$ is denoted by $\ell_j$ and the length of the block $C_j$ is denoted by $m_j$.  We denote the empty set by $\emptyset$. If $C_0=\emptyset$,    $d_1=0$ and $m_0=0$, while if $C_0\ne \emptyset$,  $d_1=1$. When $j\ge 1$, $\ell_j> 0$ and $m_j>0$. Associated with the block decomposition is the sequence 
 $m_0,\ell_1,m_1,\ell_2, m_3,\ldots$ of consecutive lengths of the blocks. This sequence leads us to introduce in $\mathbb N$ a corresponding sequence of consecutive intervals  $K_0,J_1, K_1, J_2,K_2\ldots$,  which   play an important   role.

    Put $s_0=0$.  If $C_0\ne \emptyset $ we write $C_0=d_1d_2\ldots d_{t_0}$. If $C_0=\emptyset$ we put $t_0=0$. 
    When $j\ge 1$ write 
      \begin{equation}
      B_j=d_{t_{j-1}+1}d_{t_{j-1}+2}\ldots d_{s_j}\ 
     {\rm and}\ C_j=d_{s_j+1}d_{s_j+2}\ldots d_{t_j}.\label{eq:stdefinition}
       \end{equation}
     This defines $t_0,s_0,t_1,s_1,\cdots$ and we put
     \begin{align}
     J_j&=\{t_{j-1}+1,\ldots, s_j\},  \ {\rm for}\ j\in {\mathbb N},   \  \label{eq:Jdefinition}\\
       K_0&=\emptyset, \ {\rm if}\ C_0=\emptyset,\ {\rm and} \label{eq:K0)definition}\\
       K_j&=\{s_j+1,\ldots, t_j\},  \ {\rm for}\ j\in \{0\}\cup{\mathbb N}, \ {\rm if} \ C_0\ne \emptyset. \label{eq:Kdefinition}
     \end{align}
     As $B_j$ has length $\ell_j$, $J_j$ has $\ell_j$ elements, and as $C_j$ has length $m_j$, $K_j$ has $m_j$ elements. Thus,
     \begin{equation}
     \ell_j=s_j-t_{j-1} \ {\rm for}\ j=1,2,3, \ {\rm and}\ m_j=t_j-s_j\ {\rm for}\ j=0,1,2\ldots.\label{eq:blocklengths1}
     \end{equation}
 The  sets in the family $\{J_j:j\in {\mathbb N}\}\cup \{K_k:k\in \{0\}\cup {\mathbb N}\}$ of intervals   are pairwise disjoint and their union is $\mathbb N$.    If $i\in \cup_{j=1}^{\infty}J_j$ then $d_i=0$, while if  $i\in \cup_{j=0}^{\infty}K_j$ then $d_i=1$.     Note that  $m_0=t_0$, and 
   \begin{equation}
   s_j=\sum_{u=1}^j\ell_u+\sum_{u=0}^{j-1}m_u\ {\rm for}\ j\ge 2,\  {\rm   and\  that }\   t_j=\sum_{u=1}^j\ell_u+\sum_{u=0}^jm_u\ {\rm for}\ j\ge 1.\label{eq:blocklengths2}
    \end{equation}
        When $n\in {\mathbb N}$,   there is a unique $j\in \{0\}\cup {\mathbb N}$ such that either $n\in J_j$ or $n\in K_j$. In such a case, $j$ depends on $n$, and  as $n\to\infty$ so too does $j\to \infty$, and \emph{vice versa}.

\section{Estimates for sums, averages and  orbits}

 Here we consider when $x\in \Sigma$, and we use the notations concerning the digits of $x$ introduced in the previous section.    We   consider some  preliminary lemmas.
 
  \begin{lemma} \label{lemma1}Let  $ p\in (0,\infty)$, let $t, \ell\in  {\mathbb N}$ and let $J=\{t+1,\ldots t+\ell\}$.   Let $x\in \Sigma$ and let the binary expansion of $x$ be given by
\[x=\sum_{i=1}^{\infty}\frac{d_i}{2^i}.\]
Also, assume that $d_i=0$ for all $i\in J$ and that $d_{t+\ell+1}=1$. Then if $1\le r\le \ell$,    
\begin{equation}
2^{p\ell}\le \sum_{k\in\{t+1,\ldots,t+r\}}\frac{1}{f^{k-1}(x)^p}\le \frac{2^{2p}}{2^p-1}\cdot 2^{p\ell}. \label{eq:Jinequality}
\end{equation}
   \end{lemma}
 \emph{Proof.} Let   $1\le q\le \ell$ and let $k=t+q\in J$. Then  $d_{t+1}=\cdots=d_{t+\ell}=0$ and 
 \begin{align}
 f^{k-1}(x)&=\sum_{i=1}^{\infty}\frac{d_{i+k-1}}{2^i}\nonumber\\
 &=\sum_{i=t+\ell+1}^{\infty}\frac{d_i}{2^{i-t-q+1}}\label{eq:a}\\
 &\le \sum_{i=t+\ell+1}^{\infty}\frac{1}{2^{i-t-q+1}}\nonumber\\
 &=\frac{1}{2^{\ell-q+1}}.\label{eq:b}
 \end{align}
 Also, as $d_{t+\ell+1}=1$ we have from (\ref{eq:a}) and (\ref{eq:b})
 that
 \[\frac{1}{2^{\ell-q+2}}\le f^{k-1}(x)\le \frac{1}{2^{\ell-q+1}}.\]
    Thus, if $1\le r\le \ell$, 
    \[ 2^{p\ell}
   \le\sum_{q=1}^r2^{p({\ell-q+1})}
   \le \sum_{k \in \{t+1, \ldots,t+r\}}\frac{1}{f^{k-1}(x)^p}
   \le \sum_{q=1}^r2^{p({\ell-q+2})}
       \le\frac{2^{2p}}{2^p-1} \cdot 2^{p\ell}. \]
  
   \hfill$\square$
    \begin{lemma}\label{lemma2} Let  $\in p\in (0,\infty)$, let $s\in \{0\}\cup {\mathbb N}$,  let $m\in {\mathbb N}$ and let $K=\{s+1,\ldots,s+m\}$.   Let $x\in \Sigma$ and let the binary expansion of $x$ be given by
\[x=\sum_{i=1}^{\infty}\frac{d_i}{2^i}.\]
Assume that $d_i=1$ for all $i\in K$ and that $d_{s+m+1}=0$.  Then, for all  $1\le r\le m$,  
\begin{equation}
 r\le \sum_{k\in \{s+1,\ldots, s+r\}}\frac{1}{f^{k-1}(x)^p}\le2^pr.\label{eq:Kinequality}
\end{equation}
 \end{lemma}
\emph{Proof.}  Let $1\le q\le r\le m$ be such that   $k=s+q\in K.$\      As $d_k=1$ we have      
\begin{equation}
\frac{1}{2}=\frac{d_k}{2}= \frac{d_{s+q}}{2}  \le  \sum_{i=1}^{\infty}\frac{d_{i+s+q-1}}{2^i}=f^{s+q-1}(x)  =  f^{k-1}(x)\le 
   1. \nonumber \end{equation}
  We deduce that  if $ 1\le r\le m$,  
\begin{equation}
r\le \sum_{k=s+1}^{s+r}  \frac{1}{f^{k-1}(x)^p} \le 2^pr.\nonumber
\end{equation}
\hfill$\square$

 \begin{lemma}\label{lemma3} 
 
  Let $(a_n)$ and  $(b_n)$ be two sequences  of positive numbers such that $\lim_{n\to\infty}a_n=\lim_{n\to\infty}b_n=\infty$.  Let $(v_n)$ be a sequence of positive numbers,   and let $c,d, c_1,d_1>0$ and $n_0\in {\mathbb N}$  be  such that, for all $n\in {\mathbb N}$ with $n\ge n_0$,  
 \[c_1(a_n+c)\le  v_n\le c_2(b_n+d). \]
  Then there are $c_1^{\prime},c_2^{\prime}>0$ such that, for all $n\in {\mathbb N}$,   
   \[c_1^{\prime}a_n\le v_n\le c_2^{\prime}b_n.\]

 \end{lemma}

  \emph{Proof.} It is given that for all $n\ge n_0$,  
  \begin{equation}
  c_1\left(1+\frac{c}{a_n}\right)a_n\le v_n\le c_2\left(1+\frac{d}{b_n}\right)b_n.\label{eq:inequality}
  \end{equation} 
  As $a_n\rightarrow\infty$ and $b_n\rightarrow\infty$, $ca_n^{-1}\rightarrow 0$ and $=db_n^{-1}\rightarrow 0$. Thus, it is clear that there are $c_1^{\prime}, c_2^{\prime}>0$ such that,  for all $n\in {\mathbb N}$,
   $c_1^{\prime}\le c_1(1+ca_n^{-1})$ \rm  and $c_2^{\prime}\le c_2(1+db_n^{-1})$. The result follows from (\ref{eq:inequality}).
   \hfill $\square$

  DEFINITIONS. Given $p\in (0,\infty)$ and  $x\in \Sigma$ we define associated functions   $\Phi_1, \Psi_1: {\mathbb N}\longrightarrow  (0,\infty)$.   The subsets $J_j$ and $K_j$  of $\mathbb N$ are as described by (\ref{eq:stdefinition}), (\ref{eq:Jdefinition}), (\ref{eq:K0)definition}) and (\ref{eq:Kdefinition}).   
  
  If $n\in \cup_{u=1}^{\infty}J_u$ there is a unique $j\in {\mathbb N}$ with $n\in J_j$. If $n\in J_1$,  we put   \[\Phi_1(n)= \Psi_1(n)=1,  \]
  while if  $n\in J_j$ with $j\ge 2$,   we put 
  \begin{equation}
  \Phi_1(n)= \Psi_1(n)=\sum_{u=1}^j2^{p\ell_u}+\sum_{u=1}^{j-1}m_u.\label{eq:functiona}
  \end{equation}
   
 If  $n\in \cup_{u=0}^{\infty}K_u$,  there is a unique $j\in {\mathbb N}$ with $n\in K_j$. If $n\in K_0\cup K_1$ we put 
 \[\Phi_1(n)= \Psi_1(n)=1,\]  
  while if $n\in K_j$ with $j\ge 2$ we put 
    \begin{equation}
    \Phi_1(n)=   \sum_{u=1}^j2^{p\ell_u}+\sum_{u=1}^{j-1}m_u  , \  \Psi_1(n)=\ \sum_{u=1}^j2^{p\ell_u}+\sum_{u=1}^jm_u.   \label{eq:functionb}
    \end{equation}
            
    Note that the functions $\Phi_1$ and $\Psi_1$ do not depend upon whether the initial block $C_0$ is empty or non-empty according as to whether the first digit in the binary expansion of $x$ is $0$ or $1$ respectively. Their definitions have been formulated so that, later on, the main  results do not need to be stated as separate cases.  Note also that $\Phi_1$ and $\Psi_1$ are constant on each of the intervals $J_1,J_2,\ldots$ and $K_0,K_1, K_2,\ldots$.   An inspection of the definitions reveals that $\Phi_1\le \Psi_1$. 
    
    Given $x\in \Sigma$, we will use the notations  $m_0,\ell_1,m_1,\ell_2,m_2,\ldots$,  $K_0,J_1,K_1,J_2,\ldots$, $t_0,t_1,\ldots$ and $s_0,s_1,\ldots$, as discussed in Section 2 in relation to the block decomposition of $x$, generally without further explicit reference.     
 
 \begin{theorem} \label{theorem1} Let $p\in (0,\infty)$, let $x\in \Sigma$ and let  $n\in {\mathbb N}$.   Let  $\Phi_1, \Psi_1 $ be the functions mapping $\mathbb N$ into $(0,\infty)$, as given by (\ref{eq:functiona}) and  (\ref{eq:functionb}).    Then, the following hold.
  
    (i) There are $a_1,a_2>0$ such that, for all $n\in {\mathbb N}$,  
  \begin{equation}
    a_1\Phi_1(n)  \le \sum_{k=1}^n\frac{1}{f^{k-1}(x)^p}\le a_2\Psi_1(n).\label{eq:theoremconclusion1}
    \end{equation}
    
    (ii) The following statements (a), (b) and (c)  are equivalent.
    
      \hskip 0.6cm (a)
     There is $L>0$ such that, for all $j\in {\mathbb N}$ with $j\ge 2$,
    \begin{equation}
     \frac{ m_j}{\displaystyle \sum_{u=1}^j2^{p\ell_u}+\displaystyle \sum_{u=1}^{j-1}m_u} \le L. \label{eq:firstbound}
   \end{equation}
   
   \hskip 0.6cm  (b) There is $a>0$ such that, for all $n\in {\mathbb N}$,
   \begin{equation}
   \Psi_1(n)  \le   a\Phi_1(n).\label{eq:sharpestimate0}
    \end{equation}
    
    \hskip 0.6cm 
    (c) There are $b_1,b_2>0$ such that, for all $n\in {\mathbb N}$,    \begin{equation}
    b_1\Phi_1(n)  \le \sum_{k=1}^n\frac{1}{f^{k-1}(x)^p}\le b_2\Phi_1(n).\label{eq:sharpestimate1}
    \end{equation}
    
     \hskip 0.6cm (d) There are $c_1,c_2>0$ such that, for all $n\in {\mathbb N}$,    \begin{equation}
    c_1\Psi_1(n)  \le \sum_{k=1}^n\frac{1}{f^{k-1}(x)^p}\le c_2\Psi_1(n).\label{eq:sharpestimate2}
    \end{equation}

      \end{theorem}
   
  \emph{Proof.} (i) We use Lemma \ref{lemma1} and Lemma \ref{lemma2} to  estimate the sum $\sum_{k=1}^n1/f^{k-1}(x)^p.$    
 
   \emph{Case I: $n\in J_j$.}      Using  (\ref{eq:Jdefinition}) and (\ref{eq:blocklengths1}),    as $n\in J_j$,    $n=t_{j-1}+r$ for some $1\le r\le \ell_j$. In this  case, applying the right-hand inequalities in Lemma \ref{lemma1} and Lemma \ref{lemma2} respectively to   $J_1,\ldots,J_j$ in place of $J$  and $K_0,K_1,\ldots, K_{j-1}$ in place of $K$  gives
 \begin{align}
 \sum_{k=1}^n\frac{1}{f^{k-1}(x)^p}
 &\le \sum_{k\in \cup_{u=1}^jJ_u}\frac{1}{f^{k-1}(x)^p}+\sum_{k\in \cup_{u=0}^{j-1}K_u}\frac{1}{f^{k-1}(x)^p}\nonumber\\
 &\le \frac{2^{2p}}{2^p-1} \left( \sum_{u=1}^j 2^{p\ell_u}\right)+2^p\left( \sum_{u=0}^{j-1}m_u\right)\nonumber\\
 &\le \frac{2^{2p}}{2^p-1}\left(\sum_{u=1}^j 2^{p\ell_u}+ \sum_{u=0}^{j-1}m_u\right).\label{eq:inequalityref1}
 \end{align}

 Now, observe  by (\ref{eq:Jdefinition}) that  $t_{u-1}+1\in J_u$ for $u=1,2,\ldots,j$,   so that
 \[ \Bigl\{t_0+1,t_1+1,\ldots, t_{j-1}+1\Bigr\}\cup\left(\bigcup_{u=0}^{j-1}K_u\right)\subseteq  \Bigl\{1,2,\ldots,n \Bigr\}.\]
 Consequently, as $t_{j-1}+1\le t_{j-1}+r= n$,
 \[\sum_{u=1}^j \frac{1}{f^{t_{u-1}}(x)^p}+\sum_{k\in \cup_{u=0}^{j-1}K_u}\frac{1}{f^{k-1}(x)^p}
 \le \sum_{k=1}^n\frac{1}{f^{k-1}(x)^p}.\]
 We now apply the inequality on the left hand side of  (\ref{eq:Jinequality})  in  Lemma \ref{lemma1}, successively taking $\ell_1,\ell_2,\ldots\ell_j$ in place of $\ell$,  and with $r=1$ in each case.   Also we apply the inequality on the left hand side of  (\ref{eq:Kinequality}) in  Lemma \ref{lemma2}, successively  taking  $m_0,m_1, \ldots, m_{j-1}$ in place of $m$ and   with $m_0,m_1, \ldots, m_{j-1}$ successively taken in place of $r$. We thus obtain  

 \begin{align}
 \displaystyle\sum_{u=1}^j2^{p\ell_u}+\displaystyle\sum_{u=0}^{j-1}m_u\le\sum_{k=1}^n\frac{1}{f^{k-1}(x)^p}.\label{eq:inequalityref2}
\end{align}
 Using (\ref{eq:inequalityref1}) and (\ref{eq:inequalityref2}),  we deduce that when $n\in J_j$ with $j\ge 2$,
  
  \begin{align}
 \displaystyle\sum_{u=1}^j2^{p\ell_u}+\displaystyle\sum_{u=0}^{j-1}m_u
 \le \sum_{k=1}^n\frac{1}{f^{k-1}(x)^p}\le   \frac{2^{2p}}{2^p-1}\left( \sum_{u=1}^j2^{p\ell_u}+ \sum_{u=0}^{j-1}m_u\right)\label{eq:estimate1}
  \end{align}

   Now, observe that  $j$ is a function of $n$ and that as $n\longrightarrow\infty$, $j\longrightarrow\infty$ also, in which case
    \begin{equation}
    \sum_{u=1}^j2^{p\ell_u}+ \sum_{u=1}^{j-1}m_u \longrightarrow\infty\ {\rm and}\  \sum_{u=1}^j\ell_u+\displaystyle\sum_{u=1}^{j-1}m_u\longrightarrow\infty.\label{eq:goestoinfinity}
    \end{equation}
Using (\ref{eq:estimate1}) and (\ref{eq:goestoinfinity})  and applying Lemma \ref{lemma3}  with $v_n=\sum_{k=1}^n1/f^{k-1}(x)^p$,  $a_n=b_n= \sum_{u=1}^j2^{p\ell_u}+\displaystyle\sum_{u=1}^{j-1}m_u$ and $c=d=m_0$, we deduce that  there are $a_1^{\prime},a_2^{\prime}>0$ such that for all $j\in {\mathbb N}$ with $j\ge 2$,
\[a_1^{\prime}\left(\sum_{u=1}^j2^{p\ell_u}+ \sum_{u=1}^{j-1}m_u\right)\le \sum_{k=1}^n\frac{1}{f^{k-1}(x)^p}\le a_2^{\prime} \left(\sum_{u=1}^j2^{p\ell_u}+ \sum_{u=1}^{j-1}m_u\right). \]
An inspection of the definition of $\Phi_1$ and $\Psi_1$ on $\cup_{j=1}^{\infty}J_j$ in (\ref{eq:functiona}),  shows that    \begin{equation}
  a_1^{\prime }\Phi_1(n)\le  \displaystyle\sum_{k=1}^n\frac{1}{f^{k-1}(x)^p} \le a_2^{\prime}\Psi_1(n),   \label{eq:inequA}   \end{equation}
 {\rm for all}\ $n\in \cup_{j=2}^{\infty}J_j$.
\vskip 0.2cm

  \emph{Case II: $n\in K_j.$}  Using (\ref{eq:Kdefinition}) and (\ref{eq:blocklengths1}), write $n=s_j+r$ where $1\le r\le m_j$.  We have 
 \begin{align}
& \sum_{k=1}^n\frac{1}{f^{k-1}(x)^p}\nonumber \\
&=\sum_{k\in \cup_{u=1}^jJ_u}\frac{1}{f^{k-1}(x)^p}+\sum_{k\in \cup_{u=0}^{j-1}K_u}\frac{1}{f^{k-1}(x)^p}+\sum_{k\in \{s_j+1,\ldots,s_j+r\}}\frac{1}{f^{k-1}(x)^p}.\nonumber
 \end{align}
 Using Lemma \ref{lemma1} and Lemma \ref{lemma2} on the sets $J_u, K_u$, in an  argument  similar to  Case I where  (\ref{eq:estimate1}) was obtained, 
it follows that for all  $n \in K_j$ with $j\ge 2$,
\begin{equation}
\sum_{u=1}^j2^{p\ell_u}+\sum_{j=0}^{j-1}m_u\le\sum_{k=1}^n\frac{1}{f^{k-1}(x)^p}\le\frac{2^{2p}}{2^p-1}\left(\sum_{u=1}^j2^{p\ell_u}+\sum_{u=0}^jm_u\right). \label{eq:inequB}
\end{equation}
By using Lemma \ref{lemma3}, and using the same type of argument for $K_j$ as was used to obtain (\ref{eq:inequA})  in the case of $J_j$, we deduce from (\ref{eq:inequB}) that   there are $a_1^{\prime\prime},a_2^{\prime\prime}>0$ such that
 \begin{equation}
a_1^{\prime\prime} \Phi_1(n)\le \sum_{k=1}^n\frac{1}{f^{k-1}(x)^p}\le a_2^{\prime\prime}\Psi_1(n),\label{eq:theKinequality}
 \end{equation}
 for all $n\in \cup_{j=2}^{\infty}K_j.$
 Then, (i) of Theorem \ref{theorem1} now follows from (\ref{eq:inequA}) and (\ref{eq:theKinequality}) upon noting that the definitions of $\Phi_1$ and $\Psi_1$ in (\ref{eq:functiona}) and (\ref{eq:functionb})
  mean that these functions are equal when restricted to $K_0\cup J_1\cup K_1$.
 
 \emph{Proof that (ii) (a)  implies (ii) (b).} Assume that (ii) (a)  holds.   Observe from the definitions of $\Phi_1$ and $\Psi_1$ in  (\ref{eq:functiona}) that $\Psi_1(n) /\Phi_1(n)=1$  for all $n\in \cup_{u=1}^{\infty}J_u$.  Also, if $n\in K_j$ with $j\ge 2$ we have from the definitions of $\Phi_1$ and $\Psi_1$ in  (\ref{eq:functionb}) that
 \begin{align}
   \frac{ \Psi_1(n)}{\Phi_1(n)}
  &= \frac{\displaystyle\sum_{u=1}^j2^{p\ell_u}+\displaystyle\sum_{u=1}^jm_u}{\displaystyle \sum_{u=1}^j2^{p\ell_u}+\displaystyle\sum_{u=1}^{j-1}m_u} \nonumber \\
 &=  1+ \frac{m_j}{\displaystyle\sum_{u=1}^j2^{p\ell_u} +\sum_{u=1}^{j-1}m_u}\label{eq:equivalence}\\
 &\le 1+L. \label{eq:conclusiona}
 \end{align}
Thus, as  $\Phi_1(n)=\Psi_2(n)$ for all $n\in \cup_{j=1}^{\infty}J_j$, and as (\ref{eq:conclusiona}) holds for all $n\in \cup_{j=2}^{\infty}K_j$, using Lemma \ref{lemma3}  in a like manner as before we see  that (ii) (b)  holds.

 \emph{Proof that (ii) (b) implies (ii) (a).}  Assume that (ii) (a) does not hold. Then, (\ref{eq:equivalence}) shows there is no $a>0$ such that
 $\Psi(n)\Phi_1(n))^{-1}  \le a$
 for all $n\in {\mathbb N}$. This shows that there is no $a>0$ such that  inequality (\ref{eq:sharpestimate0})  holds. Thus, if (ii) (a) does not hold, (ii) (b) does not hold. Thus, (ii) (b) implies (ii) (a).

\emph{Proof of equivalence of (ii) (a),  (ii) (b),  (ii) (c) and (ii) (d).} We have shown that (ii) (a) and (ii) (b) are equivalent. The equivalence of each of (ii) (c) and (ii) (d) with each of (ii) (a) and (ii) (b) immediately follows from  Theorem  \ref{theorem1} (i), previously proved, and from the fact that $\Phi_1\le \Psi_1$.
\newline${\ }$ \hfill$\square$

  The use of Lemma \ref{lemma3} in the preceding proof is a particular use  of the fact that   questions of asymptotic behaviour  do not depend on the first few terms.     This  fact may be used without explicit reference in some of the  later  proofs in this paper.   As well, the definitions of $\Phi_1$ and $\Psi_1$ in  (\ref{eq:functiona}) and (\ref{eq:functionb}) 
as being $1$ on $K_0\cup J_1\cup K_1$ are for convenience. The values of $\Phi_1$ and $\Psi_1$  on $K_0\cup J_1\cup K_1$ can be any positive numbers without affecting the statement of Theorem  \ref{theorem1}.  The property of $\Phi_1$ as in (\ref{eq:sharpestimate1}) says in a certain sense that when it holds, $\Phi_1$ is an `optimal' estimation of the asymptotic behaviour of the sequence   $n\longmapsto \sum_{k=1}^n1/f^{k-1}(x)^p$. This is because if $\Phi$ is any function that satisfies (\ref{eq:sharpestimate1})
 with $\Phi$ in place of $\Phi_1$,  there are $b_1,b_2>0$ such that for all $n\in {\mathbb N}$ we have
$b_1\Phi(n)\le\Phi_1(n) \le b_2\Phi(n).$ Similarly, if $\Psi$ is any function that satisfies (\ref{eq:sharpestimate2})
 with $\Psi$ in place of $\Psi_1$,  there are $c_1,c_2>0$ such that, for all $n\in {\mathbb N}$,  $c_1\Psi(n)\le\Psi_1(n) \le c_2\Psi(n).$

Note that condition (\ref{eq:firstbound}) ensures that the `optimal' estimates in (\ref{eq:sharpestimate1}) and (\ref{eq:sharpestimate2}) occur,  and is fulfilled by the  simpler condition that there is $L>0$ such that
\[m_j\left(\displaystyle\sum_{u=1}^{j-1}m_u\right)^{-1}\le L,\ {\rm for\  all}\ j\ge 2.\]
 Note that the latter condition depends only upon the numbers  $m_u$, regardless of  the numbers $\ell_u$, but it is clear that the two conditions are not equivalent.

Theorem \ref{theorem1} estimates the sum $\sum_{k=1}^{n}1/f^{k-1}(x)^p$ for $n\in {\mathbb N}$   from above and below, using the functions $\Phi_1$ and $\Psi_1$. The functions $\Phi_1$ and $\Psi_1$ are defined solely in terms of the lengths of the blocks in the block decomposition of $x\in \Sigma$. We now  derive corresponding estimates for the expression $n^{-1}\left(\sum_{k=1}^n 1/f^{k-1}(x)^p \right)$,  and consider relationships between the two sets of estimates.

DEFINITIONS. Define functions $\Phi_2, \Psi_2: {\mathbb N}\longrightarrow (0,\infty)$ as follows, using the definitions of $\Phi_1$ and $\Psi_1$ as given in (\ref{eq:functiona}) and (\ref{eq:functionb}).

If $n\in \cup_{u=1}^{\infty}J_u$, there is a unique $j\in {\mathbb N}$ with $n\in J_j$. If $n\in J_1$ put
\[\Phi_2(n)=\Psi_2(n)=1,\]
while if $n\in J_j$ with $j\ge 2$, put
\begin{equation}
\Phi_2(n)=\frac{\Phi_1(n)}{\displaystyle\sum_{u=1}^j\ell_u+\sum_{u=1}^{j-1}m_u} \ {\rm and}\ \Psi_2(n)=\frac{\Psi_1(n)}{\displaystyle\sum_{u=1}^{j-1}\ell_u+\sum_{u=1}^{j-1}m_u}.\label{eq:Jdefphisub1}
\end{equation}

If $n\in \cup_{u=0}^{\infty}K_u$, there is a unique $j\in \{0\}\cup{\mathbb N} $ with $n\in K_j$. If $n\in K_0\cup K_1$ put
\[\Phi_2(n)=\Psi_2(n)=1,\]
while if $n\in K_j$ with $j\ge 2$, put
\begin{equation}
\Phi_2(n)=\frac{\Phi_1(n)}{\displaystyle\sum_{u=1}^j\ell_u+\sum_{u=1}^jm_u}\  {\rm and}\ \Psi_2(n)=\frac{\Psi_1(n)}{\displaystyle\sum_{u=1}^{j}\ell_u+\sum_{u=1}^{j-1}m_u}.\label{eq:Kdefphisub2}
\end{equation}
Note that, as for $\Phi_1$ and $\Psi_1$,  the functions $\Phi_2$ and $\Psi_2$ do not depend  upon whether the initial block $C_0$ is empty or non-empty according as to whether the first digit in the binary expansion of $x$ is $0$ or $1$ respectively. Also, as for $\Phi_1$ and $\Psi_1$,  $\Phi_2$ and $\Psi_2$ are constant on each of the intervals $J_1,J_2,\ldots$ and $K_0, K_1, K_2,\ldots$ within $\mathbb N$.   As $\Phi_1\le \Psi_1$, an inspection of the definitions reveals that $\Phi_2\le \Psi_2$.

\begin{lemma}\label{lemma4}
Let $a,b,c>0$ with $a\ge c$. Then $a/(a+b)\ge c/(b+c)$.
\end{lemma}
\emph{Proof.}  It suffices to observe that $a/(a+b)-c/(b+c)= b(a-c)/(a+b)(b+c)\ge 0.$\hfill$\square$
 
\begin{lemma}\label{lemma5} Let $\Phi_1(n)$, $\Phi_2(n)$ be given as in  (\ref{eq:functiona}), (\ref{eq:functionb}), (\ref{eq:Jdefphisub1}) and (\ref{eq:Kdefphisub2}).  Let $j,n\in {\mathbb N}$ with $j\ge 2$.  Then     the following statements hold.

(i)   If $n$ is the minimum element of $J_j$, 
\begin{equation}
\frac{2}{m_0+3}\cdot \left(1+\frac{\ell_j}{\displaystyle\sum_{u=1}^{j-1}\ell_u+\displaystyle\sum_{u=1}^{j-1}m_u}\right)\le \frac{1}{n}\cdot \frac{\Phi_1(n)}{\Phi_2(n)}.\label{eq:ineqJphiextra}
\end{equation} 

(ii)  If   $n$ is the minimum element of $K_j$,  
\begin{equation}
\frac{3}{m_0+4}\left(1+\frac{m_j}{\displaystyle\sum_{u=1}^j\ell_u+\displaystyle\sum_{u=1}^{j-1}m_u}\right)\le \frac{1}{n}\cdot \frac{\Phi_1(n)}{\Phi_2(n)}.\label{eq:ineqKphiextra}
\end{equation}

\end{lemma}

\emph{Proof of (i).}   The minimum element of $J_j$ is $\sum_{u=1}^{j-1}\ell_u+\sum_{u=0}^{j-1}m_u+1$.  Using this value for $n$, and the definitions of   $\Phi_1(n)$  and $\Phi_2(n)$ as in  (\ref{eq:functiona}) and  (\ref{eq:Jdefphisub1}),    we see that  
        \begin{align}
\frac{1}{n} \cdot \frac{\Phi_1(n)}{\Phi_2(n)} 
      = \frac{\displaystyle\sum_{u=1}^j\ell_u+\sum_{u=1}^{j-1}m_u}{\displaystyle\sum_{u=1}^{j-1}\ell_u+\sum_{i=0}^{j-1}m_u +1 }
     =\left(1+\frac{\ell_j}{\displaystyle \sum_{u=1}^{j-1}\ell_u+\displaystyle\sum_{u=1}^{j-1}m_u}\right)\cdot \left(\frac{\displaystyle \sum_{u=1}^{j-1}\ell_u+\displaystyle\sum_{u=1}^{j-1}m_u}{\displaystyle \sum_{u=1}^{j-1}\ell_u+\displaystyle\sum_{u=0}^{j-1}m_u+1}\right).\label{eq:leftfraction2}
        \end{align}
         Putting  $a=\sum_{u=1}^{j-1}\ell_u+\sum_{u=1}^{j-1}m_u\ge 2$, $b=m_0+1$ and $c=2$,   Lemma \ref{lemma4} gives  
          \[\frac{2}{m_0+3}\le \frac{\displaystyle\sum_{u=1}^{j-1}\ell_u+\sum_{u=1}^{j-1}m_u}{ \displaystyle\sum_{u=1}^{j-1}\ell_j+\displaystyle\sum_{u=0}^{j-1}m_j+1},\]
         which together with (\ref{eq:leftfraction2})  proves (\ref{eq:ineqJphiextra}) and proves (i).
         
          \emph{Proof of (ii)}. The minimum element of $K_j$  is $\sum_{u=1}^j\ell_u+\sum_{u=0}^{j-1}m_u+1$.  Using this value for $n$ and  the definition of   $\Phi_2(n)$ as in   (\ref{eq:Kdefphisub2}),    we see that  
         \begin{align}
\frac{1}{n} \cdot \frac{\Phi_1(n)}{\Phi_2(n)} 
      = \frac{\displaystyle\sum_{u=1}^j\ell_u+\sum_{u=1}^jm_u}{\displaystyle\sum_{u=1}^j\ell_u+\sum_{u=0}^{j-1}m_u +1 }
     = \left(1+\frac{m_j}{\displaystyle \sum_{u=1}^j\ell_u+\displaystyle\sum_{u=1}^{j-1}m_u}\right) \left(\frac{\displaystyle \sum_{u=1}^j\ell_u+\displaystyle\sum_{u=1}^{j-1}m_u}{\displaystyle \sum_{u=1}^j\ell_u+\displaystyle\sum_{u=0}^{j-1}m_u+1}\right).\label{eq:leftfraction4}
        \end{align}
        Putting  $a=\sum_{u=1}^j\ell_u+\sum_{u=1}^{j-1}m_u\ge 3$,   $b=m_0+1$ and $c=3$, applying Lemma \ref{lemma4} together with ({\ref{eq:leftfraction4})   gives (\ref{eq:ineqKphiextra}) and proves (ii).

         \hfill$\square$
                
                  \begin{lemma}\label{lemma6}
          Let $\Psi_1$, $\Psi_2$ be given as in  (\ref{eq:functiona}), (\ref{eq:functionb}), (\ref{eq:Jdefphisub1}) and (\ref{eq:Kdefphisub2}).  Let $j,n\in {\mathbb N}$ with $j\ge 2$.  Then     the following statements hold.

(i) If $n$ is the maximum element of $J_j$,  
\begin{equation}
\frac{1}{n}\cdot \frac{\Psi_1(n)}{\Psi_2(n)}\le \frac{1}{1+\ell_j\left(\displaystyle\sum_{u=1}^{j-1}\ell_u+\displaystyle\sum_{u=1}^{j-1}m_u\right)^{-1}}.\label{eq:maxJ}
\end{equation}

(ii) If $n$ is the maximum element of $K_j$, 
\begin{equation}
\frac{1}{n}\cdot \frac{\Psi_1(n)}{\Psi_2(n)}\le \frac{1}{1+m_j\left(\displaystyle\sum_{u=1}^j\ell_u+\displaystyle\sum_{u=1}^{j-1}m_u\right)^{-1}}.\label{eq:maxK}
\end{equation}
\end{lemma}

\emph{Proof of (i).}   The maximum element of $J_j$ is $\sum_{u=1}^j\ell_u+\sum_{u=0}^{j-1}m_u$.  Using this value for $n$ and the definition of   $\Psi_2(n)$ as in   (\ref{eq:Jdefphisub1})  we see that   
     \begin{align*}
    \frac{1}{n} \cdot \frac{\Psi_1(n)}{\Psi_2(n)}
   =  \frac{\displaystyle\sum_{u=1}^{j-1}\ell_u+ \displaystyle\sum_{u=1}^{j-1}m_u} {\displaystyle\sum_{u=1}^j\ell_u+\displaystyle\sum_{u=0}^{j-1}m_u}
    \le\frac{1}{1+\ell_j\left(\displaystyle\sum_{u=1}^{j-1}\ell_u+\displaystyle\sum_{u=1}^{j-1}m_u\right)^{-1}} .
          \end{align*}
          \emph{Proof of (ii).} The maximum element of $K_j$ is $\sum_{u=1}^j\ell_u+\sum_{u=0}^jm_u$. Using this value for $n$ and the definition of   $\Psi_2(n)$ as in   (\ref{eq:Kdefphisub2})  gives
     \begin{align*}
    \frac{1}{n} \cdot \frac{\Psi_1(n)}{\Psi_2(n)}
   =  \frac{\displaystyle\sum_{u=1}^j\ell_u+ \displaystyle\sum_{u=1}^{j-1}m_u} {\displaystyle\sum_{u=1}^j\ell_u+\displaystyle\sum_{u=0}^jm_u}
    \le\frac{1}{1+m_j\left(\displaystyle\sum_{u=1}^j\ell_u+\displaystyle\sum_{u=1}^{j-1}m_u\right)^{-1}} .
          \end{align*}
\hfill $\square$

 \begin{lemma} \label{lemma7}
Let $p>0$ and let  $j\in {\mathbb N}$ with $j\ge 2$. Then, if $\ell_1,\ldots \ell_j, m_1,\ldots, m_j\in {\mathbb N}$,
\begin{equation}
\frac{m_j}{\displaystyle\sum_{u=1}^j2^{p\ell_u}+\displaystyle\sum_{u=1}^{j-1}m_u} \le \max\left(1,\frac{1}{p}\right)\left(\frac{m_j}{\displaystyle\sum_{u=1}^j\ell_u+\displaystyle\sum_{u=1}^{j-1}m_u}\right).\nonumber
\end{equation}
\end{lemma} 
\emph{Proof.} This follows from the fact that $2^{p\ell}\ge p\ell$.\hfill$\square$

DEFINITION. The block decomposition of a number $x\in \Sigma$ is called   \emph{regular}  if there is $C>0$ such that for all $j\ge 2$, 
\begin{equation}
\max\left(\frac{\ell_j}{\displaystyle\sum_{u=1}^{j-1}\ell_u+\displaystyle\sum_{u=1}^{j-1}m_u}, \frac{m_j}{\displaystyle\sum_{u=1}^j\ell_u+\displaystyle\sum_{u=1}^{j-1}m_u }\right)\le C.\label{eq:regularbound}
\end{equation}
The definition  restricts the possible lengths of $B_j$ and $C_j$ as $j\to\infty$, so that  they cannot   increase `too quickly' compared with the sums of the lengths of  preceding blocks, and  it expresses a notion of `balance' between the blocks.

\begin{theorem}\label{theorem2}
 Let $p\in (0,\infty)$. Let $x\in \Sigma$ and let  $n\in {\mathbb N}$.  Let $\ell_1,m_1,\ell_2,m_2,\ldots$ and $ J_1,K_1, J_2,K_2,\ldots$ and other notations  be   as  described  in Section 2. 
 Let  $\Phi_2, \Psi_2 $ be the associated functions mapping $\mathbb N$ into $(0,\infty)$, as given by (\ref{eq:Jdefphisub1}) and  (\ref{eq:Kdefphisub2}).    Then, the following hold. 
 
 (i) There are $b_1,b_2>0$ such that, for all $n\in {\mathbb N}$,
  \begin{equation}
    b_1\Phi_2(n)\le \frac{1}{n}\left( \sum_{k=1}^n\frac{1}{f^{k-1}(x)^p}\right)\le b_2\Psi_2(n).\label{eq:averagesresult}
    \end{equation}
    
    (ii) The following statements (a), (b), (c) and (d) are equivalent.

     \hskip 0.6cm (a) The block decomposition of $x$ is regular.
     
     \hskip 0.6cm (b) There is $a>0$ such that for all $n\in {\mathbb N}$,
     \begin{equation}
     \Psi_2(n)\le a\Phi_2(n). \label{eq:regularphi0}
     \end{equation}
     
      \hskip 0.6cm (c) There are $b_1,b_2>0$ such that for all $n\in {\mathbb N},$ 
    \begin{equation}
    b_1\Phi_2(n)\le \frac{1}{n}\left(\sum_{k=1}^n\frac{1}{f^{k-1}(n)^p}\right)\le b_2\Phi_2(n).\label{eq:regularphi1}
    \end{equation}
     
    \hskip 0.6cm (d) There are $c_1,c_2>0$ such that for all $n\in {\mathbb N},$ 
    \begin{equation}
    c_1\Psi_2(n)\le \frac{1}{n}\left(\sum_{k=1}^n\frac{1}{f^{k-1}(n)^p}\right)\le c_2\Psi_2(n).\label{eq:regularphi2}
    \end{equation}
    \end{theorem}
    \emph{Proof of (i).} 
      Case I.   If $n\in J_j$ with $j\ge 2$, observe from (\ref{eq:blocklengths1}) and (\ref{eq:blocklengths2}) that
   $\sum_{u=1}^{j-1}\ell_u+\sum_{u=0}^{j-1}m_u+1\le n\le \sum_{u=1}^j\ell_u+\sum_{u=0}^{j-1}m_u.$
   Then, using   Lemma \ref{lemma3}, by Theorem \ref{theorem1} (i) we deduce  that there are $c_1^{\prime},c_2^{\prime}>0$,  such that for all $n\in J_j$ with $j\ge 2$, 
       \[c_1^{\prime}\cdot \frac{\Phi_1(n)}{\displaystyle\sum_{u=1}^j\ell_u+\displaystyle\sum_{u=1}^{j-1}m_u}\le \frac{1}{n}\left( \sum_{k=1}^n\frac{1}{f^{k-1}(x)^p}\right)\le c_2^{\prime}\cdot \frac{\Psi_1(n)}{\displaystyle\sum_{u=1}^{j-1}\ell_u+\displaystyle\sum_{u=1}^{j-1}m_u}, \]
   which by again using  Lemma \ref{lemma3} and (\ref{eq:Jdefphisub1})
 shows that there are $c_1^{\prime\prime},c_2^{\prime\prime}>0$ such that for all $n\in \cup_{j=1}^{\infty}J_j$,
  \begin{equation}
  c_1^{\prime\prime}\Phi_2(n)\le  \frac{1}{n}\left( \sum_{k=1}^n\frac{1}{f^{k-1}(x)^p}\right)\le c_2^{\prime\prime}\Psi_2(n).\label{eq:XXX}
  \end{equation}
    
   Case II. If  $n\in K_j$ with $j\ge 2$, observe that $\sum_{u=1}^j\ell_u+\sum_{u=0}^{j-1}m_u+1\le n\le \sum_{u=1}^j\ell_u+\sum_{u=0}^jm_u.$
   Using  Lemma \ref{lemma3}  and   Theorem  \ref{theorem1} (i) there are $d_1^{\prime},d_2^{\prime}>0$ such that for all $n\in K_j$ with $j\ge 2$,
   \[d_1^{\prime}\cdot \frac{\Phi_1(n)}{\displaystyle\sum_{u=1}^j\ell_u+\displaystyle\sum_{u=1}^jm_u}\le\frac{1}{n}\left( \sum_{k=1}^n\frac{1}{f^{k-1}(x)^p}\right)\le d_2^{\prime}\cdot \frac{\Psi_1(n)}{\displaystyle\sum_{u=1}^j\ell_u+\displaystyle\sum_{u=1}^{j-1}m_u}, \]
   which  by again using Lemma \ref{lemma3}  and (\ref{eq:Kdefphisub2}) shows that for some $d_1^{\prime\prime}>0, d_2^{\prime\prime}>0$, for all $n\in\cup_{j=0}^{\infty} K_j$ we have
       \begin{equation}   
   d_1^{\prime\prime} \Phi_2(n)\le  \frac{1}{n}\left( \sum_{k=1}^n\frac{1}{f^{k-1}(x)^p}\right)\le d_2^{\prime\prime} \Psi_2(n).\label{eq:YYY}
   \end{equation}
   We deduce from (\ref{eq:XXX}) and (\ref{eq:YYY}) that for suitable $b_1,b_2>0$, (\ref{eq:averagesresult}) holds for all $n\in {\mathbb N}$, which proves (i).
 
   \emph{Proof that (ii) (a) and (ii) (b) are equivalent.}  A calculation using (\ref{eq:functiona})  and (\ref{eq:Jdefphisub1}) shows that 
    if  $n\in J_j$ and $j\ge 2$,  
\begin{align}
\frac{\Psi_2(n)}{\Phi_2(n)}=1+\frac{\ell_j}{\displaystyle\sum_{u=1}^{j-1}\ell_u+\sum_{u=1}^{j-1}m_u}.\label{eq:ratiopsiphi1}
\end{align}
Another     calculation using (\ref{eq:functionb}) and and (\ref{eq:Kdefphisub2})  shows that if  $n\in K_j$ and $j\ge 2$, we have
\begin{align}
 \frac{\Psi_2(n)}{\Phi_2(n)}=\left(1+\frac{m_j}{\displaystyle \sum_{u=1}^j2^{p\ell_u}+\sum_{u=1}^{j-1}m_u}\right)\cdot\left(1+\frac{m_j}{\displaystyle\sum_{u=1}^j\ell_u+\sum_{u=1}^{j-1}m_u}\right).\label{eq:ratiopsiphi2}
\end{align}
Assuming (ii) (a) holds the block decomposition is regular, so if $C$ is  as in (\ref{eq:regularbound}),  we see from (\ref{eq:ratiopsiphi1}), (\ref{eq:ratiopsiphi2})
 and Lemma \ref{lemma7}  that 
\begin{equation}
\Psi_2(n)\le \max\left(2,1+\frac{1}{p}\right)(C+1)^2\Phi_2(n),\ {\rm  for\    all}\  n\in \cup_{n=2}^{\infty}(J_j\cup K_u).\label{eq:inequalityX}
\end{equation}
 It now follows from  Lemma \ref{lemma3} and (\ref{eq:inequalityX})   that there is $a>0$ such that    
 \begin{equation}
  \Psi_2(n)\le a\Phi_2(n),\label{eq:inequalityY}
 \end{equation}
  for all $n\in {\mathbb N}$.     Thus, (ii) (a) implies (ii) (b). Conversely, if (ii) (b) holds, $\Psi_2/\Phi_2$ is bounded over $\mathbb N$, and from  (\ref{eq:ratiopsiphi1}) and (\ref{eq:ratiopsiphi2}) it follows that the block decomposition is regular. That is, (ii) (b) holds. Thus, (ii) (a) and (ii) (b) are equivalent.

\emph{Proof that (ii) (b) implies (ii) (c).}  If (ii) (b) holds,  the inequality (\ref{eq:regularphi0}) when applied to the right hand side of (\ref{eq:averagesresult}), already proved, gives (ii) (c). 

\emph{Proof that (ii) (a) implies (ii) (d).}
Inequality  (\ref{eq:inequalityY}) was derived from assumption (ii) (a). But  by writing
(\ref{eq:inequalityY}) in the form
$a^{-1}\Psi_2(n)\le \Phi_2(n)$ 
for all $n\in {\mathbb N}$,  we see that (ii) (d) follows from the left hand side of (\ref{eq:averagesresult}) in (i), already proved.

\emph{Proof that (ii) (c) implies (ii) (a).}  Let $j\ge 2$ and let $n$ be the minimum element in $J_j$.  Then if $a_1$ is the constant in (\ref{eq:theoremconclusion1}) in Theorem  \ref{theorem1} (i),  applying  (\ref{eq:ineqJphiextra}) in   Lemma \ref{lemma5} (i) gives
\begin{align}
(n\Phi_2(n))^{-1}\left(\sum_{k=1}^n\frac{1}{f^{k-1}(x)^p}\right)\ge \frac{a_1}{n}\cdot \frac{\Phi_1(n)}{\Phi_2(n)}
\ge\frac{2a_1}{m_0+3}\left(1+\frac{\ell_j}{\displaystyle\sum_{u=1}^{j-1}\ell_u+\sum_{u=1}^{j-1}m_u}\right).\label{eq:phiratio1}
\end{align} 

Similarly,  if $j\ge 2$,  $n$ is the minimum element in $K_j$, and $a_1$ is the constant in (\ref{eq:theoremconclusion1}),  Theorem  \ref{theorem1} (i) and  (\ref{eq:ineqKphiextra}) in Lemma \ref{lemma5} (ii)  give
\begin{align}
(n\Phi_2(n))^{-1}\left(\sum_{k=1}^n\frac{1}{f^{k-1}(x)^p}\right)\ge\frac{a_1}{n}\cdot \frac{\Phi_1(n)}{\Phi_2(n)}
\ge \frac{3a_1}{m_0+4}\left(1+\frac{m_j}{\displaystyle\sum_{u=1}^j\ell_u+\sum_{u=1}^{j-1}m_u}\right).\label{eq:phiratio2}
\end{align}
Then,  (\ref{eq:phiratio1}) and (\ref{eq:phiratio2}) imply  that if (ii) (a) fails, there is no $b>0$ such that  
\[\frac{1}{n} \left(\displaystyle\sum_{u=1}^{k-1}\frac{1}{f^{k-1}(x)^p}\right)\le b\Phi_2(n), {\rm for\  all}\  n.\]
Thus, if (ii) (a) fails, (ii) (c) fails. We deduce that (ii) (c) implies (ii) (a).

\emph{Proof that (ii) (d) implies (ii) (a).}  Let $j\ge 2$ and let $n$ be the maximum element in $J_j$. Then if $a_2$ is the constant in (\ref{eq:theoremconclusion1}),   Theorem \ref{theorem1} (i) and (\ref{eq:maxJ}) in   Lemma \ref{lemma6} (i)    give
  \begin{equation}
  (n\Psi_2(n))^{-1}\left(\sum_{k=1}^n\frac{1}{f^{k-1}(x)^p}\right) \le \frac{a_2}{n}\cdot\frac{\Psi_1(n)}{\Psi_2(n)}\le a_2\left(
  \frac{1}{1+\ell_j\left(\displaystyle\sum_{u=1}^{j-1}\ell_u+\displaystyle\sum_{u=1}^{j-1}m_u\right)^{-1}}\right).\label{eq:Psiratio1}
  \end{equation}
  Similarly,  if $j\ge 2$ and  $n$ is the maximum element in $K_j$,  and if $a_2$ is the constant in (\ref{eq:theoremconclusion1}),   Theorem  \ref{theorem1} (i) and  (\ref{eq:maxK}) in Lemma \ref{lemma6} (ii)  give  
  \begin{equation}
  (n\Psi_2(n))^{-1}\left(\sum_{k=1}^n\frac{1}{f^{k-1}(x)^p}\right) \le \frac{a_2}{n}\cdot\frac{\Psi_1(n)}{\Psi_2(n)}\le a_2\left(
  \frac{1}{1+m_j\left(\displaystyle\sum_{u=1}^j\ell_u+\displaystyle\sum_{u=1}^{j-1}m_u\right)^{-1}}\right).\label{eq:Psiratio2}
\end{equation}
Now,  (\ref{eq:Psiratio1}) and (\ref{eq:Psiratio2}) imply that if (ii) (a) fails, there is no $C>0$ such that  
\[\frac{1}{n} \left(\displaystyle\sum_{u=1}^{k-1}\frac{1}{f^{k-1}(x)^p}\right)\ge C\Psi_2(n), {\rm for\  all}\  n.\]
Thus, if (ii) (a) fails, (ii) (d) fails. We deduce that (ii) (d) implies (ii) (a). \hfill $\square$

       Theorem \ref{theorem1} gives one estimation on the growth of the sequence $n\longmapsto\sum_{k=1}^n1/f^{k-1}(x)^p$ and   Theorem \ref{theorem2} provides another. Specifically, the conclusions (\ref{eq:theoremconclusion1})  and (\ref{eq:averagesresult})  of Theorems \ref{theorem1} and   \ref{theorem2} respectively may be written as
 \begin{equation}
    a_1\Phi_1(n)  \le \sum_{k=1}^n\frac{1}{f^{k-1}(x)^p}\le a_2\Psi_1(n), \ {\rm and} \label{eq:comparisona}
    \end{equation}
    \begin{equation}
    b_1n\Phi_2(n)\le \sum_{k=1}^n\frac{1}{f^{k-1}(x)^p}\le b_2n\Psi_2(n),\label{eq:comparisonb}
    \end{equation}
 for all $n\in {\mathbb N}$.   Both (\ref{eq:comparisona}) and (\ref{eq:comparisonb}) provide upper and lower estimates on   $ \sum_{k=1}^n1/f^{k-1}(x)^p$ as $n\longrightarrow \infty$.   However, whereas the  estimates $\Phi_1(n)$ and $\Psi_1(n)$ are defined entirely in terms of the lengths of the blocks and the  interval $J_j$ or $K_j$ in which $n$ occurs,   the  estimates $n\Phi_2(n)$ and $n\Psi_2(n)$ are defined in terms of the lengths of the blocks and  $n$ itself.  This difference means that the estimates in (\ref{eq:comparisona}) and (\ref{eq:comparisonb}) for $  \sum_{k=1}^n1/f^{k-1}(x)^p$  are not necessarily `equivalent'.  In particular, observe that the left hand sides of (\ref{eq:comparisona}) and (\ref{eq:comparisonb}) provide two lower bounds for the function $n\longmapsto \sum_{k=1}^n1/f^{k-1}(x)^p$, given by the functions $n\longmapsto a_1\Phi_1(n)$ and $n\longmapsto b_1n\Phi_2(n)$. We might regard these lower bounds as `equivalent' if there are $c_1,c_2>0$ such that 
 \[c_1\Phi_1(n)\le n\Phi_2(n)\le c_2\Phi_1(n),\] 
 for all $n\in {\mathbb N}$.   It turns out that this occurs precisely when the block decomposition of $x$ is regular, and  the same result holds for the corresponding right hand estimates of $   \sum_{k=1}^n1/f^{k-1}(x)^p$ in  (\ref{eq:comparisona}) and (\ref{eq:comparisonb}). We need some preliminary lemmas.

\begin{lemma}\label{lemma8} Let $\Phi_1(n)$, $\Phi_2(n)$ be given as in  (\ref{eq:functiona}), (\ref{eq:functionb}), (\ref{eq:Jdefphisub1}) and (\ref{eq:Kdefphisub2}).  Let $j,n\in {\mathbb N}$ with $j\ge 2$.  Then     the following hold.

(i) If   $n\in J_j$,  
\begin{equation}
  \frac{3}{m_0+3}\le \frac{1}{n}\cdot \frac{\Phi_1(n)}{\Phi_2(n)}\le 1+\frac{\ell_j}{\displaystyle\sum_{u=1}^{j-1}\ell_u+\displaystyle\sum_{u=1}^{j-1}m_u}. \label{eq:ineqJphi}
\end{equation}

(ii) If   $n\in K_j$, 
\begin{equation}
  \frac{4}{m_0+4}\le \frac{1}{n}\cdot \frac{\Phi_1(n)}{\Phi_2(n)}\le 1+\frac{m_j}{\displaystyle\sum_{u=1}^{j-1}\ell_u+\displaystyle\sum_{u=1}^{j-1}m_u}. \label{eq:ineqKphi}
\end{equation}

\end{lemma}
\emph{Proof of (i).}    When $n\in J_j$ with $j\ge 2$, 
 $\sum_{u=1}^{j-1}\ell_u+\sum_{u=0}^{j-1}m_j+1\le n\le \sum_{u=1}^j\ell_j+\sum_{u=0}^{j-1}m_u.$
 Using  the definition of   $\Phi_2(n)$ as in   (\ref{eq:Jdefphisub1})  we see that 
\begin{align}
\frac{\displaystyle\sum_{u=1}^j\ell_u+\sum_{u=1}^{j-1}m_u}{ \displaystyle\sum_{u=1}^j\ell_u+\displaystyle\sum_{u=0}^{j-1}m_u}\le \frac{1}{n} \cdot \frac{\Phi_1(n)}{\Phi_2(n)} 
      \le \frac{\displaystyle\sum_{u=1}^j\ell_u+\sum_{u=1}^{j-1}m_u}{\displaystyle\sum_{u=1}^{j-1}\ell_u+\sum_{u=1}^{j-1}m_u  }
     =1+\frac{\ell_j}{\displaystyle \sum_{u=1}^{j-1}\ell_u+\displaystyle\sum_{u=1}^{j-1}m_u}.\label{eq:leftfraction1}
        \end{align}
        Putting  $a=\sum_{u=1}^j\ell_u+\sum_{u=1}^{j-1}m_u\ge 3$, $b=m_0$ and $c=3$, we may apply Lemma \ref{lemma4} to deduce that  
          \[\frac{3}{m_0+3}\le \frac{\displaystyle\sum_{u=1}^j\ell_u+\sum_{u=1}^{j-1}m_u}{ \displaystyle\sum_{u=1}^j\ell_u+\displaystyle\sum_{u=0}^{j-1}m_u},\]
         which together with (\ref{eq:leftfraction1}) gives (\ref{eq:ineqJphi}) and proves (i).
         
         \emph{Proof of (ii).}     When $n\in K_j$ with $j\ge 2$, 
 $\sum_{u=1}^j\ell_u+\sum_{u=0}^{j-1}m_u+1\le n\le \sum_{u=1}^j\ell_u+\sum_{u=0}^jm_u.$
  Using  the definition of   $\Phi_2(n)$ as in   (\ref{eq:Kdefphisub2})  gives
\begin{align}
\frac{\displaystyle\sum_{u=1}^j\ell_u+\sum_{u=1}^jm_u}{ \displaystyle\sum_{u=1}^j\ell_u+\displaystyle\sum_{u=0}^jm_u}\le \frac{1}{n} \cdot \frac{\Phi_1(n)}{\Phi_2(n)} 
      \le \frac{\displaystyle\sum_{u=1}^j\ell_u+\sum_{u=1}^jm_u}{\displaystyle\sum_{u=1}^j\ell_u+\sum_{u=1}^{j-1}m_u  }
     =1+\frac{m_j}{\displaystyle \sum_{u=1}^j\ell_u+\displaystyle\sum_{u=1}^{j-1}m_u}.\label{eq:leftfraction3}
        \end{align}
        As $\sum_{u=1}^j\ell_u+\sum_{u=1}^jm_u\ge4$ for $j\ge 2$, apply Lemma \ref{lemma4} with $a=\sum_{u=1}^j\ell_u+\sum_{u=1}^{j-1}m_u$, $b=m_0$ and $c=4$. This gives
        \[\frac{4}{m_0+4}\le \frac{\displaystyle\sum_{u=1}^j\ell_u+\sum_{u=1}^jm_u}{ \displaystyle\sum_{u=1}^j\ell_u+\displaystyle\sum_{u=0}^jm_u},\]
        which together with (\ref{eq:leftfraction3}) proves   (\ref{eq:ineqKphi}) and (ii). 
          \hfill$\square$

         \begin{lemma}\label{lemma9}
          Let $\Psi_1$, $\Psi_2$ be given as in  (\ref{eq:functiona}), (\ref{eq:functionb}), (\ref{eq:Jdefphisub1}) and (\ref{eq:Kdefphisub2}).  Let $j,n\in {\mathbb N}$ with $j\ge 2$.  Then the following statements hold.

(i) If   $n\in J_j$,    
\begin{equation}
\frac{3}{m_0+3}\left(\frac{1}{1+\ell_j\left(\displaystyle\sum_{u=1}^{j-1}\ell_u+\displaystyle\sum_{u=1}^{j-1}m_u\right)^{-1}}\right)  \le \frac{1}{n}\cdot \frac{\Psi_1(n)}{\Psi_2(n)}\le 1. \label{eq:Jstatement}
\end{equation} 

(ii)  If   $n\in K_j$,    
\begin{equation}
\frac{4}{m_0+4}\left(\frac{1}{1+m_j\left(\displaystyle\sum_{u=1}^j\ell_u+\displaystyle\sum_{u=1}^{j-1}m_u\right)^{-1}}\right) 
 \le \frac{1}{n}\cdot \frac{\Psi_1(n)}{\Psi_2(n)}\le 1. \label{eq:Kstatement}
\end{equation} 
         \end{lemma}
         
         \emph{Proof of (i).}  As $n\in J_j$, $\sum_{u=1}^{j-1}\ell_u+\sum_{u=0}^{j-1}m_u+1\le n\le \sum_{u=1}^j\ell_u+\sum_{u=0}^{j-1}m_u.$ Hence, putting $a=\sum_{u=1}^j\ell_u+\sum_{u=1}^{j-1}m_u\ge 3$, $b=m_0$ and $c=3$, and applying Lemma \ref{lemma4} gives
 \begin{align}
 \frac{3}{m_0+3}\left(\frac{1}{1+\ell_j\left(\displaystyle\sum_{u=1}^{j-1}\ell_u+\displaystyle\sum_{u=1}^{j-1}m_u\right)^{-1}}\right)   &\le\left(\frac{\displaystyle\sum_{u=1}^j\ell_u+ \displaystyle\sum_{u=1}^{j-1}m_u} {\displaystyle\sum_{u=1}^j\ell_u+\displaystyle\sum_{u=0}^{j-1}m_u} \right)\cdot \left(\frac{\displaystyle\sum_{u=1}^{j-1}\ell_u+\displaystyle\sum_{u=1}^{j-1}m_u}{\displaystyle\sum_{u=1}^j\ell_u+\displaystyle\sum_{u=1}^{j-1}m_u }\right)\nonumber\\
 &=\frac{\displaystyle\sum_{u=1}^{j-1}\ell_u+ \displaystyle\sum_{u=1}^{j-1}m_u} {\displaystyle\sum_{u=1}^j\ell_u+\displaystyle\sum_{u=0}^{j-1}m_u}\nonumber \\
 &\le\frac{1}{n} \cdot \frac{\Psi_1(n)}{\Psi_2(n)}, \label{eq:Jinequality1}\\
  &\le  \frac{\displaystyle\sum_{u=1}^{j-1}\ell_u+ \displaystyle\sum_{u=1}^{j-1}m_u} {\displaystyle\sum_{u=1}^{j-1}\ell_u+\displaystyle\sum_{u=0}^{j-1}m_u+1}\label{eq:Jinequality2}\nonumber\\
  &\le1. \end{align}
 Then, (\ref{eq:Jinequality1}) and (\ref{eq:Jinequality2}) show that (\ref{eq:Jstatement}) holds and (i) is proved
 
 \emph{Proof of (ii).}  As $n\in K_j$ we have $\sum_{u=1}^j\ell_u+\sum_{u=0}^{j-1}m_u+1\le n\le \sum_{u=1}^j\ell_u+\sum_{u=0}^jm_u.$ Hence, putting $a=\sum_{u=1}^j\ell_u+\sum_{u=1}^jm_u\ge 4$, $b=m_0$ and $c=4$, and applying Lemma \ref{lemma4} gives
 \begin{align}
 \frac{4}{m_0+4}\left(\frac{1}{1+m_j\left(\displaystyle\sum_{u=1}^j\ell_u+\displaystyle\sum_{u=1}^{j-1}m_u\right)^{-1}}\right)   &\le\left(\frac{\displaystyle\sum_{u=1}^j\ell_u+ \displaystyle\sum_{u=1}^jm_u} {\displaystyle\sum_{u=1}^j\ell_u+\displaystyle\sum_{u=0}^jm_u} \right)\cdot \left(\frac{\displaystyle\sum_{u=1}^j\ell_u+\displaystyle\sum_{u=1}^{j-1}m_u}{\displaystyle\sum_{u=1}^j\ell_u+\displaystyle\sum_{u=1}^jm_u }\right)\nonumber\\
 &=\frac{\displaystyle\sum_{u=1}^j\ell_u+ \displaystyle\sum_{u=1}^{j-1}m_u} {\displaystyle\sum_{u=1}^j\ell_u+\displaystyle\sum_{u=0}^jm_u}\nonumber \\
 &\le\frac{1}{n} \cdot \frac{\Psi_1(n)}{\Psi_2(n)},\label{eq:Kinequality1} \\
  &\le  \frac{\displaystyle\sum_{u=1}^j\ell_u+ \displaystyle\sum_{u=1}^{j-1}m_u} {\displaystyle\sum_{u=1}^j\ell_u+\displaystyle\sum_{u=0}^{j-1}m_u+1}\nonumber\\
  &\le1.\label{eq:Kinequality2}
  \end{align}
  Now, (\ref{eq:Kinequality1}) and (\ref{eq:Kinequality2}) show that (\ref{eq:Kstatement}) holds and (i) is proved
 \hfill$\square$

 DEFINITION. Given two functions $f,g:{\mathbb N}\longrightarrow (0,\infty)$,   $f$ and $g$ are called  \emph{comparable} if there are $a,b>0$ such that $af(n)\le g(n)\le bf(n)$ for all $n\in {\mathbb N}$. 
 
 The following result shows  that   the lower bounds  for $\sum_{k=1}^n1/f^{k-1}(x)^p$ given in   (\ref{eq:comparisona}) and (\ref{eq:comparisonb}) are comparable if and only if the block decomposition of $x$ is regular.   Similarly,  the   upper bounds for $\sum_{k=1}^n1/f^{k-1}(x)^p$ given in   (\ref{eq:comparisona}) and (\ref{eq:comparisonb}) are likewise comparable if and only if the block decomposition of $x$ is regular.
  \begin{theorem}\label{theorem3}
 Let $p\in (0,\infty)$ and let $x\in \Sigma$. The following statements (i) to (v) are equivalent.
 
 (i) The block decomposition of $x$ is regular.

 (ii) There is $a>0$ such that $\Phi_1(n)\le an\Phi_2(n)$, for all $n\in {\mathbb N}$.
 
 (iii) There is $b>0$ such that $n\Psi_2(n)\le b\Psi_1(n)$, for all $n\in {\mathbb N}$.
 
 (iv)  There are $a_1,a_2>0$ such that for all $n\in {\mathbb N}$,
  \begin{equation}
   a_1n\Phi_2(n)\le \Phi_1(n)\le a_2n\Phi_2(n).\label{eq:phinequation}  
  \end{equation}
  
 (v) There are $b_1,b_2>0$ such that for all $n\in {\mathbb N}$,
  \begin{equation}
   b_1n\Psi_2(n)\le \Psi_1(n)\le b_2 n\Psi_2(n).\label{eq:psinequation}   
 \end{equation}

 \end{theorem}

\emph{Proof.}   Assume that the  block decomposition of $x$ is regular. Then, it follows from (\ref{eq:ineqJphi}) and (\ref{eq:ineqKphi}) of Lemma \ref{lemma8} that  (\ref{eq:phinequation}) holds. Similarly,   It follows from (\ref{eq:Jstatement}) and (\ref{eq:Kstatement}) in Lemma \ref{lemma9}   that (\ref{eq:psinequation}) holds.   Thus, (i) implies (iv) and (v).

 It is immediate that (iv) implies (ii) and that (v) implies (iii). 
 
  Now, assume that (ii) holds.  Then it follows from (\ref{eq:ineqJphiextra}) and (\ref{eq:ineqKphiextra})  in Lemma \ref{lemma5} that (i) holds. Thus, (ii) implies (i). 
 
Now, assume that (iii) holds. Then it follows from (\ref{eq:maxJ}) and (\ref{eq:maxK})  in Lemma \ref{lemma6} that (i) holds. Thus, (iii) implies (i). 

The above implications are sufficient to deduce that  (i) to (v) are equivalent.
\hfill $\square$

 Note that if  the block decomposition is regular, Lemma \ref{lemma7} shows that (\ref{eq:firstbound}) holds, so  from the equivalence of (\ref{eq:firstbound}), 
 (\ref{eq:sharpestimate1}) and (\ref{eq:sharpestimate2}) in Theorem  \ref{theorem1} (ii) we deduce that $\Phi_1$ and $\Psi_1$ are comparable. However, it is clearly possible for the block decomposition of $x$  to  satisfy (\ref{eq:firstbound}) without it being regular, but in this case Theorem \ref{theorem1} shows that  $\Phi_1$ and $\Psi_1$ are still comparable.  Theorem \ref{theorem1} shows that    $\Phi_1$ and $\Psi_1$ are  comparable if and only if inequality (\ref{eq:firstbound}) holds.

  \section{Normal and simply normal numbers}

 Let $x\in\Sigma$, and let the sequence of digits in the binary expansion of $x$ be $d_1d_2d_3\ldots$. Given $n\in {\mathbb N}$ and  a block $e_1e_2\cdots e_r$ of $r$ digits,   make the definition that
 \begin{equation}
 A( e_1e_2\ldots e_r,n)=\Bigl\{k:1\le k\le n\ {\rm and}\ d_k=e_1, d_{k+1}=e_2, \ldots,d_{k+r-1}=e_r\Bigr\}. \nonumber
 \end{equation}
 In particular,
  \begin{align}
  A(0, n)&=\Bigl\{k:1\le k\le n\ {\rm and}\ d_k= 0\Bigr\},\ {\rm and}\nonumber\\
 A(1, n)&=\Bigl\{k:1\le k\le n\ {\rm and}\ d_k= 1\Bigr\}. \nonumber
  \end{align}
 We will use the notation that if $A$ is a finite set, $|A|$ denotes its number of elements.  Then $x$ is called  a \emph{normal number} if for all $r$   and any block $ e_1e_2\cdots e_r$ of length $r$,   
 \begin{equation}
 \lim_{n\to\infty}\frac{1}{n}\cdot \left|A( e_1e_2\ldots e_r,n)\right|=\frac{1}{2^r}.\label{eq:normalnumber}
 \end{equation}
 That is, $x$ is normal if, for each $r$ and each block $E$  of length $r$, the asymptotic proportion of occurrences of $E$ in the expansion of $x$ is the same as the probability of obtaining $E$ as the outcome of selecting at random $r$ zeros  or ones. Strictly speaking, the above defines a \emph{normal number to the base} $2$}, but as we do not consider other bases, we simply use the term \emph{normal number}.
 
 If (\ref{eq:normalnumber}) holds for $r=0$ and $r=1$, $x$ is called \emph{simply normal}. Thus,  $x$ is simply normal precisely when
 \begin{equation}
  \lim_{n\to\infty}\frac{1}{n}\cdot \left|A(0,n)\right|=\frac{1}{2}\ {\rm and}\ \lim_{n\to\infty}\frac{1}{n}\cdot \left|A(1,n)\right|=\frac{1}{2}.\label{eq:simplynormal}
  \end{equation}
  As $A(0,n)+A(1,n)=n$,   each of the conditions in (\ref{eq:simplynormal}) implies the other. 
\begin{theorem} \label{theorem4}
 Let $p>1$ and let $x\in\Sigma$. Then if  $x$ is normal,
 \[\lim_{n\to\infty}\frac{1}{n}\left(\sum_{k=1}^n\frac{1}{f^{k-1}(x)^p} \right)=\infty. \]
\end{theorem}

 \emph{Proof.}  One can check  that if $r\in {\mathbb N}$ and $e_1e_2e_3\ldots$ is an infinite sequence of zeros and ones containing no infinite sequence of consecutive ones,   then  
\[ \sum_{i=1}^{\infty}\frac{e_i}{2^i}<\frac{1}{2^r} \ {\rm if\  and\  only\  if}\ e_1=e_2=\cdots=e_r=0.\]
 Let the binary expansion of $x=\sum_{i=1}^{\infty}d_i/2^i$. Then, as  $f^{k-1}(x)$ $=\sum_{i=1}^{\infty}d_{i+k-1}/2^i,$
 we see that 
 $f^{k-1}(x)\in \left[0, 1/2^r\right)\ {\rm if\  and\  only\  if}\ d_k=d_{k+1}=\cdots=d_{k+r-1}=0.$ Thus, 
 \begin{equation}
 f^{k-1}(x)\in \left[0, \frac{1}{2^r}\right) \ {\rm if\  and\  only\  if}\  k\in A(00\cdots0, n),\label{eq:normalinterval}
 \end{equation}
 where $A(00\cdots0, n)$ has $r$ zeros. As $x$ is normal, we see from (\ref{eq:normalinterval}) that
\[\lim_{n\to\infty}\frac{1}{n}\cdot\Bigg|\left\{k: 1\le k\le n\ {\rm and}\ f^{k-1}(x)\in \left[0, \frac{1}{2^r}\right) \right\} \Bigg| =\lim_{n\to\infty}\frac{1}{n}\cdot |A(00\cdots0, n)|=\frac{1}{2^r}.\] 
So,   for each $r$, we see that for all sufficiently large $n$,
\[\frac{1}{n}\cdot\Bigg|\left\{k: 1\le k\le n\ {\rm and}\ f^{k-1}(x)^p\in \left[0, \frac{1}{2^{rp}}\right) \right\} \Bigg| >\frac{1}{2^{r+1}}.\]
Now let $M>0$. As $p>1$,    choose $r$ so that  $2^{r(p-1)-1}>M$.   We deduce that,   for all sufficiently large values of  $n$,  
\[\frac{1}{n}\left(\sum_{k=1}^n\frac{1}{f^{k-1}(x)^p}\right)>\frac{2^{rp}}{2^{r+1}}=2^{r(p-1)-1}>M. \]
 \hfill$\square$
 
 A number may be simply normal but not normal. The following gives a characterisation of simply normal numbers in terms of the lengths of the blocks in the binary expansion of the number. 
  
\begin{theorem}\label{theorem5}
Let $x\in \Sigma$, and  let $\ell_1,\ell_2,\ldots$ and $m_0,m_1,m_2,\ldots$ respectively be the lengths of the blocks of zeros and ones  in the block decomposition of $x$, as described in Section 2.  Then the following conditions (i), (ii) and (iii) are equivalent.

(i) $x$ is simply normal.

(ii) \begin{equation}
 \lim_{j\to\infty}\frac{\ell_1+\cdots+\ell_j}{m_1+\cdots+m_j}=1\ 
   and\ \lim_{j\to\infty}\frac{\ell_j}{\ell_1+\cdots+\ell_{j-1}}=0.\nonumber
    \end{equation}

 (iii)\begin{equation}
  \lim_{j\to\infty}\frac{\ell_1+\cdots+\ell_j}{m_1+\cdots+m_j}=1\ 
   and\ \lim_{j\to\infty}\frac{m_j}{m_1+\cdots+m_{j-1}}=0.\nonumber
  \end{equation}
  
  Furthermore, the following statements (iv)
  and (v) hold.
  
  (iv) If $x$ is simply normal the block decomposition of $x$ is regular.

  (v) When $p>0$, when $x$ is simply normal, and  when $\Phi_2$ is given by (\ref{eq:Jdefphisub1}) and (\ref{eq:Kdefphisub2}),  there are  $d_1,d_2>0$ such that
    \begin{equation}
    d_1\Phi_2 (n)\le \frac{1}{n}\left( \sum_{k=1}^n\frac{1}{f^{k-1}(x)^p}\right)\le d_2\Phi_2(n), \  for \ all\ n\in {\mathbb N}.\label{eq:simplynormalbound}
\end{equation}
 \end{theorem}
 
\emph{Proof.} Letting (i) hold, we prove (ii) and (iii). As $x$ is simply normal,  
\begin{equation}\lim_{j\to\infty}\frac{\ell_1+\cdots+\ell_j}{\ell_1+\cdots+\ell_j+m_0+m_1+\cdots +m_j}=\frac{1}{2}, \nonumber\end{equation} 
which gives\begin{equation}
\lim_{j\to\infty}(\ell_1+\cdots+\ell_j)(m_1+\cdots+m_j)^{-1}=1. \label{eq:balance1} 
\end {equation} 
Again, as $x$ is simply normal, 
\begin{equation}
\lim_{j\to\infty}\frac{\ell_1+\cdots+\ell_{j+1}}{\ell_1+\cdots+\ell_{j+1}+m_0+m_1+\cdots+m_j}=\frac{1}{2}, \nonumber
\end{equation}
which gives
\begin{equation}
\lim_{j\to\infty}(\ell_1+\cdots+\ell_{j+1})(m_1+\cdots+m_j)^{-1}=1. \label{eq:balance2}
\end{equation}
Using (\ref{eq:balance1}) and (\ref{eq:balance2}),   
\begin{equation}
\lim_{j\to\infty}\frac{\ell_1+\cdots+\ell_j+\ell_{j+1}}{\ell_1+\cdots+\ell_j}=1\ {\rm and\ so}\  
\lim_{j\to\infty}\frac{\ell_j}{\ell_1+\cdots+\ell_{j-1}}=0. \label{eq:ratio1}
\end{equation}  
Now (\ref{eq:balance1}) and (\ref{eq:ratio1}) show that (ii) holds so that (i) implies (ii).  

Now, still assuming that (i) holds,  put $j-1$ in place of $j$ in  (\ref{eq:balance2}) to obtain
\begin{equation}
\lim_{j\to\infty}(\ell_1+\cdots+\ell_j)(m_1+\cdots+m_{j-1})^{-1}=1.\label{eq:ratio2}
\end{equation}
Using (\ref{eq:balance1}) and (\ref{eq:ratio2}) shows that
\begin{equation}
\lim_{j\to\infty}\frac{m_1+\cdots+m_j}{m_1+\cdots+m_{j-1}}=1\ {\rm and}\  
\lim_{j\to\infty}\frac{m_j}{m_1+\cdots+m_{j-1}}=0. \label{eq:ratio3}
\end{equation}
Then (\ref{eq:balance1}) and (\ref{eq:ratio3}) show that (iii) holds so that (i) implies (iii). Thus, (i) implies both (ii) and (iii).

\emph{Proof that (ii) implies (i)}.  
 Firstly,  consider when  $n\in J_j$ with $j\ge 2$, in which case $\sum_{u=1}^{j-1}\ell_u+   \sum_{u=0}^{j-1}m_u+1\le n\le \sum_{u=1}^j\ell_u+   \sum_{u=0}^{j-1}m_u$. Also, $\sum_{u=1}^{j-1}\ell_u  \le A(0,n)\le\sum_{u=1}^j\ell_u$.
Hence,
\begin{equation}
 \left(\sum_{u=1}^{j-1}\ell_u\right) \left(\sum_{u=1}^j\ell_u+ \sum_{u=0}^{j-1} m_u \right)^{-1} \le \frac{1}{n}\cdot A(0,n)\le
  \left(\sum_{u=1}^j\ell_u\right) \left(\sum_{u=1}^{j-1}\ell_u+ \sum_{u=1}^{j-1} m_u \right)^{-1}. \label{eq:inequ2}
 \end{equation}
Since (ii) holds,
\begin{align}
\left(
\displaystyle\sum_{u=1}^{j-1}\ell_u\right) \left(\displaystyle\sum_{u=1}^j\ell_u+\displaystyle\sum_{u=0}^{j-1} m_u\right) ^{-1}
  &= \left(1+\ell_j\left(\displaystyle\sum_{u=1}^{j-1}\ell_u \right)^{-1}+ \left(\displaystyle\sum_{u=1}^{j-1}\ell_u\right)^{-1}\left(\displaystyle\sum_{u=0}^{j-1}m_u
  \right)\right)^{-1}
  \nonumber \\
  &\longrightarrow  \frac{1}{2} \ {\rm as}\ j\longrightarrow\infty.\label{eq:inequ3}
 \end{align}
   Again, as (ii) holds,   
    \begin{align}
   \left(\displaystyle\sum_{u=1}^j\ell_u\right)\left(\displaystyle\sum_{u=1}^{j-1}\ell_u+\displaystyle\sum_{u=1}^{j-1} m_u\right)^{-1}
 &= \left(1+\ell_j\left(\displaystyle\sum_{u=1}^{j-1}\ell_u\right)^{-1}\right) \left(1+\left(\displaystyle\sum_{u=1}^{j-1}\ell_u\right)^{-1}\left(\displaystyle\sum_{u=1}^{j-1}m_u\right)\right)^{-1}\nonumber\\
 &\longrightarrow  \frac{1}{2} \ {\rm as}\ j\longrightarrow\infty.\label{eq:inequ4}
 \end{align}
   
 Now, if we consider  (\ref{eq:inequ2}), (\ref{eq:inequ3}) and (\ref{eq:inequ4}) as $n\to\infty$ through values in $\cup_{j=2}^{\infty}J_j$, we see that $ n^{-1}A(0,n)\to 1/2$.
 
 Secondly, consider when $n\in K_j$ with $j\ge 2$.  In this case, we have $ \sum_{u=1}^j\ell_u +\sum_{u=0}^{j-1}m_u+1\le n\le \sum_{u=1}^j\ell_u +\sum_{u=0}^jm_u$ and $A(0,n)=\sum_{u=1}^j\ell_u$. 
 Thus,  
  \begin{equation}
\left(\sum_{u=1}^j\ell_u \right)\left(\sum_{u=1}^j\ell_u+ \sum_{u=0}^j m_u \right)^{-1}  \le \frac{1}{n}\cdot A(0,n)\le\left(\sum_{u=1}^j\ell_u\right) \left(\sum_{u=1}^j\ell_u+\displaystyle\sum_{u=1}^{j-1} m_u \right)^{-1}. \label{eq:inequ6}
 \end{equation}
     As (ii) holds,   
  \begin{equation}
   \left(\sum_{u=1}^j\ell_u\right) \left(\sum_{u=1}^j\ell_u+\sum_{u=0}^j m_u \right)^{-1}\longrightarrow \frac{1}{2}\ {\rm  as}\  j\longrightarrow\infty.\label{eq:inequ7}
   \end{equation}
   Again, as (ii) holds,
   \begin{align}
   \left(\sum_{u=1}^j\ell_u\right)& \left(\sum_{u=1}^j\ell_u+\sum_{u=1}^{j-1}m_u \right)^{-1}\nonumber\\
   &=\left(1+ \left(1+\ell_j\left(\sum_{u=1}^{j-1}\ell_u\right)^{-1}\right)^{-1}\left(\sum_{u=1}^{j-1}\ell_u\right)^{-1} \left(\sum_{u=1}^{j-1}m_u\right) \right)^{-1}\nonumber\\
   &\longrightarrow \frac{1}{2}, \ {\rm as}\  j\to\infty.\label{eq:inequ8}
   \end{align}
    If we consider  (\ref{eq:inequ6}), (\ref{eq:inequ7}) and  (\ref{eq:inequ8}) as $n\to\infty$ through values in $\cup_{j=2}^{\infty}K_j$, we see that $  n^{-1}A(0,n)\to 1/2$. 
    
    Now, as  it has been  shown that $ n^{-1}A(0,n)\to 1/2$ as $n\to\infty$ through values in either $\cup_{j=1}^{\infty}J_j$ or $\cup_{j=0}^{\infty}K_j$, it follows that $A(0,n)\to 1/2$ as $n\to\infty$, and that $x$ is simply normal. Thus, that (ii) implies (i).
    
    \emph{Proof that (iii) implies (i)}.  Consider when $j\ge 2$ and  $n\in J_j$. Then, $\sum_{u=1}^{j-1}\ell_u+\sum_{u=0}^{j-1}m_u+1\le n\le \sum_{u=1}^j\ell_u+\sum_{u=0}^{j-1}m_u$ and  $A(1,n)=\sum_{u=0}^{j-1}m_u$. Corresponding to (\ref{eq:inequ2})
 we have, for $n\in J_j$,
     \begin{align}
     \left(\sum_{u=0}^{j-1}m_u\right)\left(\sum_{u=1}^j\ell_u+ \sum_{u=0}^{j-1} m_u\right)^{-1}&\le n^{-1}A(1,n)\nonumber\\
     &\le \left(\sum_{u=0}^{j-1}m_u\right)\left(\sum_{u=1}^{j-1}\ell_u+ \sum_{u=0}^{j-1} m_u\right)^{-1}.\label{eq:inequ9}
      \end{align}
     Now, as (iii) holds,
          \begin{align}
     & \left(\sum_{u=0}^{j-1}m_u\right)\left(\sum_{u=1}^j\ell_u+\displaystyle\sum_{u=0}^{j-1} m_u\right)^{-1}\nonumber\\
       &=\left(\left(\sum_{u=0}^jm_u\right)^{-1}\left(\sum_{u=1}^j\ell_u\right) \left(1+ m_j\left(\sum_{u=0}^{j-1}m_u\right)^{-1}\right)+1\right)^{-1}\nonumber\\
     &\longrightarrow 1/2, \ {\rm as}\  j\longrightarrow \infty.\label{eq:limitA}
     \end{align}
      Again, as (iii) holds, 
      \begin{align}
      \left(\sum_{u=0}^{j-1}m_u\right)\left(\sum_{u=1}^{j-1}\ell_u+ \sum_{u=0}^{j-1} m_u\right)^{-1}&=\left(\left(\sum_{u=0}^{j-1}m_u\right)^{-1} \left(\sum_{u=1}^{j-1}\ell_u\right)+1\right)^{-1}\nonumber\\
     & \longrightarrow 1/2, \ {\rm as}\ j\longrightarrow \infty. \label{eq:limitB}
      \end{align}\
       It now follows from (\ref{eq:inequ9}), (\ref{eq:limitA}) and (\ref{eq:limitB}) that $n^{-1}A(1,n)\longrightarrow 1/2$ as $n\longrightarrow \infty$ through values in $\cup_{j=1}^{\infty}J_j.$

     Consider now when $j\ge 2$ and $n\in K_j$. Then, $\sum_{u=1}^j\ell_u+\sum_{u=0}^{j-1}m_u+1\le n\le \sum_{u=1}^j\ell_u+\sum_{u=0}^jm_u$ and  $\sum_{u=0}^{j-1}m_u\le A(1,n)\le\sum_{u=0}^jm_u$. Then,  for $n\in K_j$, 
          \begin{equation}
     \left(\sum_{u=0}^{j-1}m_u\right)\left(\sum_{u=1}^j\ell_u+ \sum_{u=0}^j m_u\right)^{-1}\le n^{-1}A(1,n)\le
     \left(\sum_{u=0}^jm_u\right)\left(\sum_{u=1}^j\ell_u+ \sum_{u=0}^{j-1} m_u\right)^{-1}.\label{eq:inequ10}
     \end{equation}
         Now, as (iii) holds,
     \begin{align}
      &\left(\sum_{u=0}^{j-1}m_u\right)\left(\sum_{u=1}^j\ell_u+ \sum_{u=0}^j m_u\right)^{-1}\nonumber\\
      &= \left[
      \left(\sum_{u=0}^jm_u\right)^{-1}
      \left(\sum_{u=1}^j\ell_u\right) \left(1+m_j\left(\sum_{u=0}^{j-1}m_u\right)^{-1}\right)
      +1+m_j\left(\sum_{u=0}^{j-1}m_u\right)^{-1}
      \right]^{-1},\nonumber\\
      &\longrightarrow1/2, \ {\rm as}\ j\longrightarrow \infty.\label{eq:tending1}
     \end{align}
       Again, as (iii) holds,
     \begin{align}
      &\left(\sum_{u=0}^jm_u\right)\left(\sum_{u=1}^j\ell_u+ \sum_{u=0}^{j-1} m_u\right)^{-1}\nonumber \\
      &=\left[\left(\sum_{u=0}^jm_u\right)^{-1}\left(\sum_{u=1}^j\ell_u\right)+\left(1+m_j\left(\sum_{u=0}^{j-1}m_u\right)^{-1}\right)^{-1} \right]^{-1},\nonumber\\
      &\longrightarrow 1/2, \ {\rm as}\ j\longrightarrow \infty.\label{eq:tending2}
     \end{align}
     It now follows from (\ref{eq:inequ10}), (\ref{eq:tending1}) and (\ref{eq:tending2}) that $n^{-1}A(1,n)\longrightarrow 1/2$ as $n\longrightarrow \infty$ through values in $\cup_{j=0}^{\infty}K_j.$
     
     Since $n^{-1}A(0,n)\longrightarrow 1/2$ as $n\to\infty$ through values  all of which  belong  either to  $\cup_{j=1}^{\infty}J_j$ or  to $\cup_{j=0}^{\infty}K_j$, it follows that $\lim_{n\to\infty}n^{-1}A(1,n)=1/2$, and so $x$ is simply normal. That is, (iii) implies (i).

              \emph{Proof of (iv) and (v).}  If 
         $x$ is simply normal, (ii) and (iii) hold and it follows from (\ref{eq:regularbound}) that the block decomposition of $x$ is regular and (iv) holds.  So, by (ii) (a) implies (ii) (c) in Theorem \ref{theorem2},  (\ref{eq:simplynormalbound}) holds and so (v) is proved.          
   ${\ }$  \hfill $\square$

 \begin{lemma}\label{lemma10}  Let $x$ be simply normal.  Then there are $c,d>0$ such that for all $j,n\in {\mathbb N}$  with $n\in J_j \cup K_j$,  
 \begin{equation}
 c\le \frac{n}{\ell_i+\ell_2+\cdots+\ell_j}\le d.\label{eq:lambda0}
 \end{equation}
   \end{lemma}
      \emph{Proof.}  Let $j\ge 2$. If $n\in J_j\cup K_j$,  $\sum_{u=1}^{j-1}\ell_{j-1}+\sum_{u=0}^{j-1}m_u+1\le n\le \sum_{u=1}^j\ell_u+\sum_{u=0}^jm_u$. Then,
      \begin{equation}
      \frac{1+\left(\displaystyle\sum_{u=1}^{j-1}\ell_u\right)^{-1}\left(\displaystyle\sum_{u=0}^{j-1}m_u\right) }{1+\ell_j\left(\displaystyle\sum_{u=1}^{j-1}\ell_u\right)^{-1}}=\frac{\displaystyle\sum_{u=1}^{j-1}\ell_u+\displaystyle\sum_{u=0}^{j-1}m_u}{\displaystyle\sum_{u=1}^j\ell_u}\le\frac{n}{\displaystyle\sum_{u=1}^j\ell_u}\le1+ \frac{\displaystyle\sum_{u=0}^jm_u}{\displaystyle\sum_{u=1}^j\ell_u}.\label{eq:expression}
      \end{equation}
     As $x$ is simply normal, by Theorem \ref{theorem5} (ii) the  left hand and  right hand expressions in (\ref{eq:expression}) converge as $n,j\to\infty$  with $n\in J_j\cup K_j$, and each has limit $2$.  Then,  noting  (\ref{eq:expression}) again,  it follows that (\ref{eq:lambda0}) holds and the lemma is proved. \hfill$\square$
      
      DEFINITIONS. Let $p\in (0,\infty)$, let $x\in \Sigma$ and  let $J_1,J_2,\ldots$, $K_0,K_1,\ldots$, $\ell_1,\ell_2,\ldots$ and $m_0,m_1,m_2,\ldots$ be  as described in Section 2.      Functions $\Lambda_1: {\mathbb N}\rightarrow (0,\infty)$  and $\Lambda_2: {\mathbb N}\rightarrow (0,\infty)$  are defined by putting  $\Lambda_1(n)=\Lambda_2(n)=1$ when $n\in K_0$, while if  $n\ge 1$ and $n\in J_j\cup K_j$    put
\begin{equation}
\Lambda_1(n)= \sum_{u=1}^j2^{p\ell_u}\ {\rm and}\  \Lambda_2(n)= \frac{\displaystyle\sum_{u=1}^j2^{p\ell_u}}{\displaystyle\sum_{u=1}^j\ell_u} .\label{eq:lambda}
\end{equation}
      
       \begin{theorem}\label{theorem6}
Let  $p\in (0,\infty)$, let $x\in \Sigma$ be simply normal and let $\Lambda_1$, $\Lambda_2$ be the functions on $\mathbb N$ as  given in   (\ref{eq:lambda}).  Then there are $c_1,c_2, d_1,d_2>0$ such that, for  all $n\in {\mathbb N}$,
 \begin{equation}c_1\Lambda_1(n)\le \sum_{k=1}^n\frac{1}{f^{k-1}(x)^p}\le c_2\Lambda_1(n), \ {\rm and}\ \label{eq:conclusion01}
    \end{equation} 
 \begin{equation}d_1\Lambda_2(n)\le \frac{1}{n}\left(\sum_{k=1}^n\frac{1}{f^{k-1}(x)^p}\right)\le d_2\Lambda_2(n).\label{eq:conclusion02}
    \end{equation} 
    \end{theorem}
    \emph{Proof.} Note that   $2^{p\ell}\ge p\ell$ for all $\ell\in {\mathbb N}$. Thus, from (\ref{eq:lambda}) we have $\Lambda_1(n)\ge p\bigl(\sum_{u=1}^j\ell_u\bigr)$ for all $n\in\cup_{j=1}^{\infty} J_j\cup K_j$. We  consider the case when $n\in\cup_{j=2}^{\infty}J_j\cup K_j$.   There is a unique $j$ with $n\in J_j\cup K_j$. Note that $j$ depends on $n$ and that $j\to\infty$ as $n\to\infty$.         Letting $\Phi_1(n)$ be as in (\ref{eq:functiona}) and (\ref{eq:functionb}), we see that
      \begin{align}
 \Lambda_1(n)  =\sum_{u=1}^j2^{p\ell_u}
 \le \Phi_1(n)=
  \left(\sum_{u=1}^j2^{p\ell_u}\right)\,
  \left(1+\frac{\displaystyle\sum_{u=1}^{j-1}m_u}{\displaystyle\sum_{u=1}^j2^{p\ell_u}}\right)
  \le\Lambda_1(n)\left(
  1+\frac{1}{p} \cdot 
  \frac{\displaystyle\sum_{u=1}^{j-1}m_u}{\displaystyle\sum_{u=1}^j\ell_u}\label{eq:wedge}
  \right). 
  \end{align}
        As $x$ is simply normal, from  Theorem \ref{theorem5} (ii) we have
   \begin{equation}
\lim_{j\to\infty}\frac{\displaystyle\sum_{u=1}^{j-1}m_u}{\displaystyle\sum_{u=1}^j\ell_u} =\lim_{j\to\infty}\left(\frac{\displaystyle\sum_{u=1}^{j-1}m_u}{\displaystyle\sum_{u=1}^{j-1}\ell_u}\right) \cdot  \left(\frac{1}{1+\ell_j\left(\displaystyle\sum_{u=1}^{j-1}\ell_u\right)^{-1}}\right)=1.\label{eq:LMlimit}
  \end{equation}
   Hence,  from (\ref{eq:wedge})  and (\ref{eq:LMlimit})   we see that  $\Lambda_1$ and $\Phi_1$ are comparable, and (\ref{eq:conclusion01}) follows from Theorem \ref{theorem1}.

    Now, to prove (\ref{eq:conclusion02}),   observe from the definitions that when $n\in J_j\cup K_j$, $\Lambda_1(n)/\Lambda_2(n)=\sum_{u=1}^j\ell_u$. Then, it is an immediate deduction from Lemma \ref{lemma10} that $\Lambda_1$ and $n\longmapsto n\Lambda_2(n)$   are comparable.       The conclusion  (\ref{eq:conclusion02}) now follows from  (\ref{eq:conclusion01}), already proved. 
     \hfill$\square$  
 
DEFINITION.  Let $p>1$ be given. Then, relative to $p$, define  $\Gamma$ to be the set of all   numbers $x$  in $\Sigma$    for which, depending upon $x$, there exists $c>0$   such that for all $n\in {\mathbb N}$, 
   \begin{equation}
   \frac{1}{n} \left(\displaystyle \sum_{k=1}^n\frac{1}{f^{k-1}(x)^p}\right)\le c.\label{eq:boundedbyn}
   \end{equation}
   
    \begin{theorem}  \label{theorem7} Let $p>1$ be given.    Then, no normal number belongs to $\Gamma$,  $\Gamma$ has measure zero, and the set of   simply normal numbers that are in $\Gamma$ form  an uncountable set.            
   \end{theorem}
   \emph{Proof.} If $x$ is  normal, Theorem \ref{theorem4} shows that  (\ref{eq:boundedbyn}) cannot hold. Thus, if $x$ is normal, $x\notin \Gamma$. However, the set of normal numbers in $\Sigma$ has Lebesgue measure $1$ (see \cite{Ni1}  and  \cite[p. 69]{Ni} for proofs without use of  measure theory, and also see \cite[p. 79]{Ku}, for example). Thus $\Gamma$, being a subset of the complement in $\Sigma$ of the set of normal numbers, must have measure zero.

Now  let  $\Gamma^{\prime}$ be the set of all numbers $x\in \Sigma$  that have a block decomposition $B_1C_1B_2C_2\ldots $  where  the lengths of all blocks are bounded by a constant depending on $x$, and where the length of block $B_j$ equals the length of block $C_j$ for all $j$.    If $x\in \Gamma^{\prime}$, if $x\in J_j\cup K_j$ with $j\ge 2$, and  if $M$ is a bound for the lengths of the blocks,   we see  that
\[\lim_{j\to\infty}\ell_j\left(\displaystyle\sum_{u=1}^{j-1}\ell_u\right)^{-1}\le\lim_{j\to\infty}\frac{M}{ j-1} =0.\]  
Thus, by (ii) of Theorem \ref{theorem5}, $x$ is simply normal. Also,  with $\Lambda_2$  given as  in (\ref{eq:lambda}), it follows  that   
\[\Lambda_2(n)\le j2^{pM} \left(\displaystyle
\sum_{u-1}^{j-1}\ell_u\right)^{-1}\le j(j-1)^{-1} 2^{pM}\le 2^{pM+1},\]
 for all $n\in {\mathbb N}$,  
and from (\ref{eq:conclusion02}) in Theorem \ref{theorem6}  it follows that $x\in \Gamma^{\prime}$.   That $\Gamma^{\prime}$ is uncountable follows from the observation that there is an uncountable number of possible block decompositions $B_1C_1B_2C_2\ldots $ where the lengths of the  blocks are consecutively equal  in pairs and  have a common bound, and from the fact that  the numbers in $\Sigma$ are uniquely determined by their block decompositions.
 \hfill $\square$ 
 
{\bf Example.} Here is an example where Theorem \ref{theorem6} is used to estimate explicitly  the asymptotic behaviour of $\sum_{k=1}^n1/f^{k-1}(x)^p$. Let $p>0$ and consider $x\in \Sigma$ where the lengths of the blocks  in the block decomposition of $x$ are $1,1, 2,2,3,3,\ldots$, and where the first block consists of zeros. Thus, $\ell_j=m_j=j$ for all $j\in {\mathbb N}$.  Observing that 
\[\lim_{j\to\infty}\frac{\ell_1+\cdots+\ell_j}{m_1+\cdots+m_j}=1\ {\rm and}\ \lim_{j\to\infty}\frac{\ell_j}{\ell_1+\ell_2+\cdots+\ell_{j-1}}=\lim_{j\to\infty}\frac{2j}{j(j-1)}=0,\]
  Theorem \ref{theorem5} (ii)  shows  that $x$ is simply normal. Now, with $\Lambda_1$ as given by (\ref{eq:lambda}),  we have that for $n\in J_j\cup K_j$,
\[\Lambda_1(n)=\sum_{u=1}^j2^{pu} =\frac{2^p}{2^p-1}(2^{pj}-1).\]
Hence,   there are  $C,D>0$, independent of $j$ and $n$ but depending on $p$, such that
for $n\in J_j\cup K_j$,
\begin{equation}
C \,2^{pj}\le \Lambda_1(n)\le D \,2^{pj}.\label{eq:lambdaexample}
\end{equation}

 Recall that $j$  depends on $n$ and is such that $n\in J_j\cup K_j$.  Now,
 $J_j\cup K_j=\{(j-1)j+1,\ldots, j(j+1)\}. $ 
 Hence,  
\[j(j-1)+1\le n\le j(j+1),\]
 which gives
\begin{equation}
 \sqrt{n}-1< j< \sqrt{n}+1.\label{eq:squarerootestimate}
 \end{equation}
  Now,   from (\ref{eq:lambdaexample}) and (\ref{eq:squarerootestimate})  
 it is seen  that for all $n,j\in{\mathbb N}$ with $n\in J_j\cup K_j$,
   \begin{align}
2^{-p}C\, 2^{p{\sqrt n}}\le C\, 2^{pj}\le \Lambda_1(n)
  \le D\, 2^{pj}\le  2^pD\,2^{p{\sqrt n}}. \label{eq:example1}
\end{align}
Thus, by (\ref{eq:conclusion01}) of Theorem \ref{theorem6}  and by  (\ref{eq:example1}),   there are $E,F>0$ such that
\[E\,2^{p{\sqrt n}}\le\sum_{k=1}^n\frac{1}{f^{k-1}(x)}\le F\,2^{p{\sqrt n}},\]
for all $n\in {\mathbb N}$. Consequently,   $x\notin \Gamma$.

\section{Conclusion}

  When $p>0$ and  $x\in \Sigma$,  the behaviour of the sum $\sum_{k=1}^n1/f^{k-1}(x)^p$ as $n\to\infty$  involves a notion of  `balance'  between the consecutive blocks of zeros and ones in the binary expansion of $x$.     A block of zeros makes a `geometric' or exponential contribution to the sum, while a block of ones makes  makes an `arithmetic' contribution.  The  totality of these contributions depends upon the lengths  of the blocks in the block decomposition of $x$. In considering $\sum_{k=1}^n1/f^{k-1}(x)^p$, if it is to be smaller to offset the exponential effects of the blocks of zeros, there  must be longer blocks of ones -- this is how consecutive points in the orbit of $f$ can  `resist moving towards  zero'.  The functions $\Phi_1$ and $\Psi_1$ in Theorem \ref{theorem1} estimate in  (\ref{eq:theoremconclusion1}) the values of $\sum_{k=1}^n1/f^{k-1}(x)^p$ for all $n$.    It should be noted that $\Phi_1$  generally provides a superior lower estimate  compared with the `obvious' lower estimate in (\ref{eq:asharpcase}).  The condition (\ref{eq:firstbound})  in Theorem \ref{theorem1}  specifies precisely when $\Psi_1$ may be replaced by  $\Phi_1$ in (\ref{eq:theoremconclusion1}), and identifies when $\Phi_1$ is effectively `equivalent' to $\Psi_1$ in the upper estimate of  $\sum_{k=1}^n1/f^{k-1}(x)^p$.    Theorem \ref{theorem2} carries out results corresponding  to Theorem \ref{theorem1}, but for estimating $n^{-1}\bigl(\sum_{k=1}^n1/f^{k-1}(x)^p\bigr)$ instead of $\sum_{k=1}^n1/f^{k-1}(x)^p$.  In considering estimates for $n^{-1}\bigl(\sum_{k=1}^n1/f^{k-1}(x)^p\bigr)$,  a  major role is played by the  notion  in (\ref{eq:regularbound}) of a regular block decomposition,   whose definition has a symmetry and again contains a notion of `balance' between the blocks.   
  
  When $x$  is simply normal to the base $2$,  the  functions $\Phi_1$ and $\Psi_1$ are  comparable and provide `equivalent' or `optimal' estimates for upper and lower bounds on $\sum_{k=1}^n1/f^{k-1}(x)^p$ as $n\to\infty$. Also, with $x$  simply normal, $\Phi_2$ and $\Psi_2$ are  comparable and provide `optimal' estimates for upper and lower bounds on $n^{-1}\bigl(\sum_{k=1}^n1/f^{k-1}(x)^p\bigr)$.   When $x$  is simply normal, for purposes of upper and lower estimates in Theorems \ref{theorem1} and \ref{theorem2}, the  functions $\Psi_1$ and $\Psi_2$ can be replaced by $\Phi_1$ and $\Phi_2$ respectively. In this case the estimates take a simplified form, as expressed by Theorem \ref{theorem6} in terms of the functions $\Lambda_1$ and $\Lambda_2$.  
  
  Hardy and Littlewood \cite[p. 251]{HL} comment that  results such as (\ref{eq:HL}) may apply more generally for certain  irrational numbers that have  continued fractions expansions with bounded coefficients, and such numbers form a  set of measure zero  \cite[p. 69]{Ki}. The result in Theorem \ref{theorem7}, that the set $\Gamma$ has measure zero, is complementary.

  \vskip 1cm  
    \noindent Rodney Nillsen
  
  \noindent School of Mathematics and Applied Statistics
  
  \noindent University of Wollongong
  
  \noindent Northfields Avenue
  
  \noindent Wollongong NSW 2522
  
   \noindent  Australia
   
   \noindent  Email: nillsen@uow.edu.au
   
   \noindent  Website: \url{https://documents.uow.edu.au/~nillsen/}


\begin{thebibliography}{9}
\bibitem{HL}G. Hardy and J. E. Littlewood, Some Problems of Diophantine Approximation: A Series of Cosecants, \textit{Bulletin of the Calcutta Mathematical Society}, \textbf{ XX} (1930), pp.251-266.
\bibitem{Ki} A. Ya. Khintchine,  \textit{Continued Fractions}, translated by Peter Wynn, Noordhoff, Groningen, The Netherlends, 1963.
\bibitem{Ku}L. Kuipers and H. Niederreiter, \textit{Uniform Distribution of Sequences}, Wiley Interscience, 1974.
\bibitem{Ni1}R. Nillsen, \textit{Normal Numbers Without Measure Theory},    Amer.  Math.  Monthly, \textbf{107} (2000), pp. 639-644.
\bibitem{Ni}R. Nillsen, \textit{Randomness and Recurrence in Dynamical Systems}, Carus Mathematical Monograhs \textbf{31}, Mathematical Association of America, Washington DC, 2010.
  \bibitem{St}K. Stromberg, \textit{An Introduction to Classical Real Analysis}, Wadsworth, Belmont, 1981.  
 \bibitem{Wa}D. D. Wall, \textit{On normal numbers}, Ph. D. thesis, Univ. of California, Berkeley, 1949. 
  \end{thebibliography}
\end{document}